\documentclass[moor,nonblindrev]{informs1} 

\usepackage{amsmath, amssymb, amsbsy, enumitem,xcolor, xfrac, amsfonts, booktabs, stackrel}
\usepackage[utf8]{inputenc}
\usepackage[margin=3cm]{geometry}
\usepackage{charter}
\usepackage{tikz}
\usepackage{multirow}

\newcommand{\rank}{\operatorname{rank}}
\newcommand{\Span}{\operatorname{span}}
\newcommand{\conv}{\operatorname{conv}}
\newcommand{\supp}{\operatorname{supp}}

\newcommand\floor[1]{\lfloor#1\rfloor}

\newcommand{\bA}{{A}}
\newcommand{\bB}{{B}}
\newcommand{\bC}{{C}}
\newcommand{\bD}{{D}}
\newcommand{\bE}{{E}}
\newcommand{\bF}{{F}}
\newcommand{\bM}{{M}}
\newcommand{\bS}{{S}}
\newcommand{\bW}{{W}}
\newcommand{\bO}{{O}}

\newcommand{\R}{\mathbb{R}}
\newcommand{\Z}{\mathbb{Z}}
\newcommand{\BE}{\mathbf{e}}
\newcommand{\ba}{\mathbf{a}}
\newcommand{\bb}{\mathbf{b}}
\newcommand{\bc}{\mathbf{c}}
\newcommand{\cc}{\mathfrak{c}}
\newcommand{\cg}{\mathfrak{g}}
\newcommand{\cp}{\mathfrak{p}}
\newcommand{\bd}{\mathbf{d}}

\newcommand*\bl{\ensuremath{\boldsymbol\ell}}
\newcommand{\cM}{\mathcal{M}}
\newcommand{\bu}{\mathbf{u}}
\newcommand{\bw}{\mathbf{w}}
\newcommand{\bv}{\mathbf{v}}
\newcommand{\bx}{\mathbf{x}}
\newcommand{\by}{\mathbf{y}}
\newcommand{\bz}{\mathbf{z}}
\newcommand{\bgamma}{\boldsymbol{\gamma}}
\newcommand{\bdelta}{\boldsymbol{\delta}}
\newcommand{\bomega}{\boldsymbol{\omega}}
\newcommand{\bbeta}{\boldsymbol{\beta}}
\newcommand{\bzero}{\mathbf{0}}
\newcommand{\bone}{\mathbf{1}}

\newcommand{\IP}{\operatorname{IP}}
\newcommand{\LP}{\operatorname{LP}}

\usepackage{natbib}
 \NatBibNumeric
 \def\bibfont{\small}%
 \bibpunct[, ]{[}{]}{,}{n}{}{,}%

\usepackage[colorlinks=true,breaklinks=true,bookmarks=true,urlcolor=blue,
     citecolor=blue,linkcolor=blue,bookmarksopen=false,draft=false]{hyperref}
     
\usepackage{footnotebackref}
     
\def\EMAIL#1{\href{mailto:#1}{#1}}

\TheoremsNumberedThrough     

\EquationsNumberedThrough    

\begin{document}


\RUNAUTHOR{Lee et al.}

\RUNTITLE{Polynomial bounds on the column number}
\TITLE{Polynomial upper bounds on the number of differing columns of $\Delta$-modular integer programs}

\ARTICLEAUTHORS{%
\AUTHOR{Jon Lee}
\AFF{Department of Industrial and Operations Engineering, University of Michigan, USA, \EMAIL{jonxlee@umich.edu}}
\AUTHOR{Joseph Paat}
\AFF{Sauder School of Business, University of British Columbia, BC Canada, \EMAIL{joseph.paat@sauder.ubc.ca}}
\AUTHOR{Ingo Stallknecht}
\AFF{Department of Mathematics, ETH Z\"{u}rich, Switzerland, \EMAIL{ingo.stallknecht@ifor.math.ethz.ch}}
\AUTHOR{Luze Xu}
\AFF{Department of Industrial and Operations Engineering, University of Michigan, USA, \EMAIL{xuluze@umich.edu}}
} 

\ABSTRACT{We study integer-valued matrices with bounded determinants. 
Such matrices appear in the theory of integer programs (IP) with bounded determinants.
%
For example, Artmann et al. showed that an IP can be solved in strongly polynomial time if the constraint matrix is bimodular, that is, the determinants are bounded in absolute value by two.
Determinants are also used to bound the $\ell_1$-distance between IP solutions and solutions of its linear relaxation. 
One of the first works to quantify the complexity of IPs with bounded determinants was that of Heller, who identified the maximum number of differing columns in a totally unimodular matrix.
Each extension of Heller's bound to general determinants has been super-polynomial in the determinants or the number of equations. 
We provide the first column bound that is polynomial in both values. 
For integer programs with box constraints, our result gives the first $\ell_1$-distance bound that is polynomial in the determinants and the number of equations.
Our result can also be used to derive a bound on the height of Graver basis elements that is polynomial in the determinants and the number of equations.
Furthermore, we show a tight bound on the number of differing columns in a bimodular matrix; this is the first tight bound since Heller.
Our analysis reveals combinatorial properties of bimodular IPs that may be of independent interest.
}


\maketitle

\section{Introduction.}

The feasible region of an integer linear program with box constraints can be written as
\[
\IP := \{\bx \in \Z^n \colon \bA\bx = \bb, ~ \bl \le \bx \le \bu\},
\]
for a constraint matrix $\bA \in \Z^{m\times n}$ with $\rank \bA = m$ and vectors $\bb \in \Z^m$ and $\bl , \bu \in (\Z \cup \{\pm \infty\})^n$ with $\bl < \bu$. 
Integer programs have been used for many decades to model problems in operations research, computer science, and mathematics; see~\cite{CCZ2014,NWbook,AS1986} and the references therein.
One parameter that impacts the structure of $\IP$ is the largest absolute $m\times m$ minor of $\bA$, which we denote by
\[
\Delta(\bA) := \max\left\{|\det \bB| \colon \bB ~\text{is an}~m\times m~\text{submatrix of}~\bA\right\}.
\]
We say that $\bA$  is {\bf$\Delta$-modular} if $\Delta(\bA) \le \Delta$.

To illustrate the impact that $\Delta(\bA)$ has on $\IP$, consider the distance between $\IP$ and its linear relaxation 
\[
\LP := \{\bx \in \R^n \colon \bA\bx = \bb, ~  \bl \le \bx \le \bu\}.
\]
This distance, which is referred to as {\bf proximity} in the literature, is defined as the maximum distance from any vertex of $\LP$ to the closest feasible $\IP$ solution:
\[
\pi := \max_{\substack{\bx^*~\text{is a vertex}\\\text{of}\ \LP}}~~ \min_{\bz^* \in\IP} ~~\|\bx^* - \bz^*\|_1.
\]
We assume $\IP \neq \emptyset$ whenever discussing $\pi$.
Proximity is often used in the analysis of integer programming algorithms. 
For instance, proximity can also be used to bound the state space of a dynamic program~\cite{EW2018}.
Proximity also translates into an upper bound on the integrality gap: for an objective vector $\bc \in \R^n$ and vectors $\bx^* \in \LP$ and $\bz^* \in \IP$ that maximize $\bx \to \bc^\intercal \bx$ over~$\LP$ and~$\IP$, respectively, the integrality gap $|\bc^\intercal \bx^* - \bc^\intercal \bz^*|$ is at most $\pi\cdot \|\bc\|_{\infty}$.
In a seminal paper by Cook et al.~\cite{CGST1986}, they showed that $\pi \le n^2 \Delta(\bA)$.\footnote{Cook et al.'s original result considers inequality-form polyhedra, the $\ell_{\infty}$ rather than $\ell_1$-distance, and \emph{totally $\Delta$-modular matrices}, which have all absolute minors bounded by $\Delta$. A closer analysis revealed that $\Delta(\bA)$ suffices;  see~\cite[Lemma 3]{LPSX2020}.}
Eisenbrand and Weismantel~\cite{EW2018} proved $\pi \le m(2m\|\bA\|_{\infty}+1)^m$, where $\|\bA\|_{\infty}$ is the largest absolute entry of $\bA$; this was the first upper bound on $\pi$ that was independent of $n$.
Their proof approach extends\footnote{See the footnote on Page 3 of~\cite{OPW2020} for a discussion on this extension.} 
to show 
\[
\pi \le m(2m+1)^m\Delta(\bA).
\]
In the special case when $\bl = \bzero$ and $\bu \equiv \infty$, Lee et al.~\cite{LPSX2020} demonstrated that $\pi \le 3m^2\log_2(\sqrt{2m}\Delta(\bA)^{1/m})\Delta(\bA)$;
their proof crucially relied on sparsity results that are not applicable when $\bl$ and $\bu$ take general values~\cite{ADOO2017}.
No upper bounds on $\pi$ have been provided that are polynomial in $\Delta(\bA)$ and $m$ for general values of $\bl$ and $\bu$.

Testing if $\IP \neq \emptyset$ is NP-hard in general~\cite{C1971}, although it can be tested in polynomial time if $n$ is fixed~\cite{Lenstra1983}.
The parameter $\Delta(\bA)$ is also known to influence how efficiently we can test if $\IP \neq \emptyset$, at least when $\Delta(\bA)$ is small.
For example, every vertex of $\LP$ is integer valued when $\Delta(\bA) = 1$.
Therefore, testing if $\IP \neq \emptyset$ simplifies to testing if $\LP \neq \emptyset$.
Matrices with $\Delta(\bA)=1$ are called {\bf unimodular}, and after elementary row operations they are equivalent to totally unimodular (TU) matrices.
The study of TU matrices dates back to Hoffman and Kruskal~\cite{HK1956} with one prominent example being the vertex-edge incidence matrix of a bipartite graph.
It remains an open question if $\IP \neq \emptyset$ can be tested efficiently when $\Delta(\bA)$ is fixed.

If $\Delta(\bA) = 2$, then $\bA$ is called {\bf bimodular}.
One prominent example of a bimodular matrix is the vertex-edge incidence matrix of a graph whose so-called odd cycle packing number is one; see~\cite{CFHJW2020,CFHW2020} for combinatorial optimization algorithms over such graphs.
When the constraint matrix is bimodular, the vertices of $\LP$ may not be integer valued.
However, such matrices do impose the nice property that if $\IP \neq \emptyset$, then every vertex of $\LP$ lies on an edge containing a vector in $\IP$~\cite{VC2009}; Veselov and Chirkov used this property in a polynomial time algorithm to test if $\IP \neq \emptyset$ when $\bA$ contains no $m\times m$ minors equal to zero.
Artmann et al.~\cite{AWZ2017} used a more combinatorial approach to design an optimization algorithm that runs in strongly polynomial time for general bimodular matrices.
The algorithm in~\cite{AWZ2017} heavily relies on Seymour's combinatorial characterization of TU matrices~\cite{S1980}. 
Cevallos et al.~\cite[Theorem 5.4]{CWZ2018} argue that compact linear extended formulations (LEFs) do not always exist for bimodular integer programs; in their paper, they write
{\it ``A natural approach to solve bimodular integer programs would have been to try to find a compact LEF of the feasible solutions to (conic) bimodular integer programs, thus avoiding the partially involved combinatorial techniques used in~\cite{AWZ2017}, which is so far the only method to efficiently solve bimodular integer programs. Theorem 5.4 shows that this approach cannot succeed."}
%
Glanzer et al. used combinatorial properties of the constraint matrix to optimize over $\IP$ efficiently in the setting when $\bA$ has at most three distinct absolute determinants~\cite{abc2021}.
These examples illustrate the importance of combinatorial properties of the constraint matrix; this leads to a third, more combinatorial, property of $\IP$ that is influenced by $\Delta(\bA)$: the number of differing columns that $\bA$ can have.

We say that two vectors $\ba,\overline{\ba} \in \R^d$ {\bf differ} if $\ba \neq \pm \overline{\ba}$.\footnote{Glanzer et al.~\cite{GWZ2018} used the term \emph{distinctive} rather than differing.}
We say that a multiset $X \subseteq \R^d$ has {\bf differing columns} if every pair of vectors in $X$ differs and $\bzero \not \in X$.
We also treat the matrix $\bA$ as a multiset of its columns, and we denote the number of differing columns by $|\bA|$.
For $m,\Delta \in \Z_{\ge 1}$, we denote the maximum number of differing columns in a $\Delta$-modular matrix by
\[
\cc(\Delta, m) := \max\left\{ |A| \colon \bA \subseteq \Z^{m},\ \rank \bA = m,~\text{and}~\Delta(\bA) \le \Delta\right\}.
\]
The $\rank A = m$ condition is necessary. 
Otherwise, one can add a row of all zeros to any integer-valued matrix with $m-1$ rows; the resulting matrix $A$ will have $\Delta(A) = 0$.

One of the first bounds on $\cc(\Delta,m)$ is due to Heller~\cite{H1957}, who proved $\cc(1,m) = \sfrac{1}{2}\cdot(m^2+m)$. 
 Early generalizations of Heller's result focused on $\cc(\Delta, m)$ for fixed values of $\Delta$.
Lee showed $\cc(\Delta,m)\le \mathfrak{f}_L(\Delta) \cdot m^{2\Delta}$~\cite[Proposition 10.1]{Lee1989} for some function $\mathfrak{f}_L$, and
Anstee showed $\cc(\Delta,m)\in  O(m^{2\Delta(1+\log_2\Delta)})$ for totally $\Delta$-modular matrices~\cite[Theorem 3.2]{A1990}.
In the case when $A$ only contains primitive columns, Kung showed $|A|\le \mathfrak{f}_K(\Delta) \cdot m^2$ for a super-polynomial function $\mathfrak{f}_K$~\cite[Theorem 1.1]{K1990}, and $|\bA| \le m^2$ when no nonzero minor is divisible by three~\cite[Theorem 1.1]{KunMat1990}.
Oxley and Walsh recently showed $|A| \le \sfrac{1}{2}\cdot(m^2+m)+m-1$ when $A$ is bimodular, but only when $A$ contains primitive columns, and only when $m$ is sufficiently large~\cite[Theorem 1.1]{OX2021}.
The best known upper bound on $\cc(\Delta,m)$ for fixed $\Delta$ is given by Geelen et al.~\cite[Theorem 2.2.4]{GNW2021}.
They demonstrated that $\cc(\Delta,m) \le \sfrac{1}{2}\cdot  m^2 + \mathfrak{f}_G(\Delta) m$, where $\mathfrak{f}_G(\Delta)$ is a number that can be lower bounded by the ``power tower'' with base $\Delta$ iterated three times.

For fixed values of $m$, the best known upper bound on $\cc(\Delta, m)$ was due to Heller~\cite{H1957} and Glanzer et al.~\cite{GWZ2018}:
%
%
\begin{equation}\label{eqOldBounds}
\cc(\Delta, m) \le 
\begin{cases}
\sfrac{1}{2}\cdot(m^2+m), & \text{if}~ \Delta = 1;\\[.1 cm]
m^2  \Delta, & \text{if}~ \Delta =2,3;\\[.1 cm]
\sfrac{1}{2}\cdot m^2  \Delta^{2+\log_2\log_2(\Delta)}, &\text{if}~\Delta \ge 4.
\end{cases}
\end{equation}
The inequality $\cc(3,m) \le 3m^2$ is present in the analysis in~\cite[Subsection 3.3]{GWZ2018} but not stated.
In summary, neither Geelen et al. nor Glanzer et al. provided an upper bound on $\cc(\Delta,m)$ that is polynomial in $\Delta$.

An interesting variation of $\cc(\Delta,m)$ is considered by Oxley and Walsh~\cite{OX2021} and Kung~\cite{KunMat1990,K1990}, who
considered the maximum number of differing primitive columns in a $\Delta$-modular matrix, which we denote  by $\cc^{{\rm p}}(\Delta,m)$; a {\bf primitive} vector $\bv =(v_1, \dotsc, v_t)$ is an integer valued vector with $\gcd\{v_1, \dotsc, v_t\} = 1$.
It is easy to see that $\cc^{{\rm p}}(1,m) = \cc(1,m)$, and only when $\Delta \ge 2$ is there a distinction between $\cc^{{\rm p}}(\Delta,m)$ and $\cc(\Delta,m)$. 
By identifying excluded minors in matroids representable by bimodular matrices, Oxley and Walsh gave a bound of $\cc^{{\rm p}}(2,m) = \cc(1,m)+m-1$ for sufficiently large values of $m$.
Our analysis shows $\cc^{{\rm p}}(2,m) = \cc(1,m)+m$ for $m\in\{3,5\}$ and $\cc^{{\rm p}}(2,m) = \cc(1,m)+m-1$ otherwise; for a lower bound, see our tight example analysis in Section~\ref{secDelta=2}, and
for a matching upper bound, see~\eqref{eqe1e2}~\eqref{eqe1e2e3},~\eqref{eqe1e2e3e4} in Section~\ref{secDelta=2}.
Of course, the big open question in this line of work is the determination of $\cc(\Delta,m)$ and $\cc^{{\rm p}}(\Delta,m)$ for general values of $\Delta$.


\subsection{Statement of results.}

Our first main result is the exact value of $\cc(2,m)$.
This is the first tight column number bound since Heller's result. 
\begin{theorem}[An exact bound when $\Delta =2$]\label{thm:bimodular}
For every $m \in \Z_{\ge 1}$, we have 
\[
\cc(2,m)  = \frac{1}{2}  (m^2+m)+m.
\]
\end{theorem}

Our proof of Theorem~\ref{thm:bimodular} reveals new combinatorial properties about bimodular matrices.
We show that submatrices contain at most one disjoint {\bf circuit}, which is an inclusion-wise minimal set of linearly dependent columns; this generalizes a result of Heller that certain submatrices of TU matrices have no circuits.
See Lemma~\ref{lemStructuralCircuits} for a precise statement and Section~\ref{secCircuitStructure} for more discussion.
As previously quoted, combinatorial properties of the constraint matrix are critical in algorithms designed for $\Delta$-modular $\IP$s. 
For this reason, we believe our combinatorial analysis may be of independent interest.

Our proof of Theorem~\ref{thm:bimodular} requires a lower bound on $\cc(2,m)$.
We give a bound for general $\Delta$.
\begin{proposition}[A lower bound on $\cc(\Delta,m)$]\label{prop:lowerbounds}
For every $\Delta, m \in \Z_{\ge 1}$, we have 
\[
\frac{1}{2} (m^2 +m)+m(\Delta -1) \le \cc(\Delta ,m).
\]  
\end{proposition}

Geelen et al.'s result implies that $\cc(\Delta,m) = \cc(1,m)+\mathfrak{h}(\Delta) m$ for some function $\mathfrak{h}$. 
Heller's result and Theorem~\ref{thm:bimodular} support our conjecture that $\mathfrak{h}(\Delta)=\Delta-1$;
we prove this when $m \le 2$.
\begin{proposition}[An exact bound when $m \le 2$]\label{prop:m2}
Suppose $m \le 2$.
For every $\Delta \in\Z_{\ge 1}$, we have
\[
\cc(\Delta,m)  = \frac{1}{2} (m^2+m)+m(\Delta-1).
\]
\end{proposition}

Our second main result is the first upper bound on $\cc(\Delta, m)$ that is polynomial in $\Delta$ and $m$.

\begin{theorem}[An upper bound on $\cc(\Delta,m)$]\label{thmLuze}
    For every $\Delta,m \in \Z_{\ge 1}$, we have 
    \begin{equation}\label{eqLuzeBound}
     \cc(\Delta, m)
   ~ \le ~
    \begin{cases}
    \sfrac{1}{2}\cdot  (m^2+m)+m(\Delta-1), & \text{if}~ \Delta \le 2;\\[.1 cm]
    \sfrac{1}{2}\cdot  (m^2 +m) \Delta^2, &\text{if}~\Delta \ge 3.
    \end{cases}
    \end{equation}
\end{theorem}

Our third main results connects $\cc(\Delta,m)$ with the proximity value $\pi$.
We apply Theorem~\ref{thmLuze} to establish the first upper bound on $\pi$ that is polynomial in $m$ and $\Delta(\bA)$.
Unlike in~\cite{LPSX2020}, our new bound applies when the variable bounds $\bl$ and $\bu$ are arbitrary.
\begin{theorem}[LP to IP proximity]\label{thmProx}
Set $\Delta := \Delta(\bA)$, where $A\in \Z^{m\times n}$ is the constraint matrix in $\IP$.
The proximity value $\pi$ satisfies 
\[
\pi \le (m+1) \Delta (2\cc(\Delta,m)+1) \le m(m+1)^2 \Delta^3 +(m+1) \Delta.
\]
\end{theorem}
Column number bounds can also be applied to bound so-called Graver basis elements in test sets for integer programs; see, e.g.,~\cite[\S 3.7]{DHK2012}.
By directly substituting Theorem~\ref{thmLuze} into the results in~\cite[\S 3.7]{DHK2012}, one derives a bound of $O(m^3\Delta^3)$ on the $\ell_1$-norm of Graver basis vectors; the previously known bound of $O(m^m\Delta)$ can be found in Diaconis et al.~\cite[Theorem 1]{DGS1993} or by modifying a proof of Eisenbrand et al.~\cite[Lemma 2]{EHK2018}.

The paper proceeds as follows.
We begin with a proof of Proposition~\ref{prop:lowerbounds} because it is used to establish the equations in Theorem~\ref{thm:bimodular} and Proposition~\ref{prop:m2}; see Section~\ref{seclowerbounds}.
Our new combinatorial results for bimodular matrices are given in Section~\ref{secCircuitStructure}.
Sections~\ref{secDelta=2},~\ref{secm=2},~\ref{secGen}, and~\ref{secProx} contain the proofs of Theorem~\ref{thm:bimodular}, Proposition~\ref{prop:m2}, Theorem~\ref{thmLuze}, and Theorem~\ref{thmProx}, respectively.

\medskip
\noindent{\bf Notation and preliminaries.}
We use bold font to denote vectors in dimension two or higher.
$\bzero$, $\bone$, and $\mathbb{I}_k$ denote the all-zero matrix, the all-one matrix, and the $k\times k$ identity matrix, respectively. 
We denote the $i$th standard unit vector in $\R^t$ by $\BE^i_t$. 
For the $i$th standard unit vector in $\R^m$, 
we drop the subscript and write $\BE^i$.
We use $\ba ^\intercal$ to denote a row vector.
We write $ [\bb^1|\cdots|\bb^t]$ to denote a matrix with columns $\bb^1, \dotsc,\bb^{t-1}$, and $\bb^t$. 
We often partition the rows of a matrix, and it is convenient to refer to these inline; for matrices $\bB \in \Z^{r\times t}$ and $\bC \in \Z^{s\times t}$, we adopt the notation
\[
(\bB, \bC) := \left[\begin{array}{l}\bB\\\bC\end{array} \right] \in \Z^{(r+s)\times t}.
\]
A {\bf basis} is an invertible (square) matrix.
%
%
We let $\conv \bB$ and $\Span \bB$ denote the convex hull of and the linear space spanned by the columns of $\bB \subseteq \R^d$, respectively.
For $\bv = (v_1,\dotsc, v_t) \in \R^t$, we denote the support of $\bv$ by $\supp \bv := \{i = 1, \dotsc, t\colon v_i \neq 0\}$.
A $\Delta$-modular matrix $\bB$ with differing columns is {\bf maximal} if there does not exist a $\Delta$-modular matrix $\bB' \supsetneq \bB$ with differing columns.

We use {\bf elementary operations} to refer to elementary row operations that preserve integrality.
Elementary operations do not affect differing columns or $\Delta$-modularity of a matrix.
%
We write $\bB \sim \bB'$ if $\bB$ and $\bB'$ are equivalent up to elementary operations.
We also freely swap columns and multiply them by $-1$ because these operations do not affect differing columns or $\Delta$-modularity.

We often analyze the determinant structure of matrices with linearly dependent rows. 
To do this, we note that every matrix $\bB \in \Z^{m\times n}$ can be transformed via elementary operations into a matrix $(\overline{\bB}, \bzero)$, where $\overline{\bB} \in \Z^{\rank \bB \times n}$ has full row rank. 
We always use $\overline{\bB}$ to denote a full row rank projection of $\bB$ obtained via elementary operations.
Elementary operations preserve linear relationships, so the following holds:
\begin{equation}\label{eqProject}
\text{Suppose that $[\bB|\bb] \sim \left[
\begin{array}{c@{\hskip.1cm}|@{\hskip.1cm}c}
\overline{\bB} & \overline{\bb}\\
\bzero & \bzero
\end{array}\right]$ and let $\bv \in \R^{|\bB|}$.
We have $\bB\bv = \bb$ if and only if $\overline{\bB} \bv = \overline{\bb}$.}
\end{equation}
Consequently,
\begin{equation}\label{eqPreserve}
\text{if $\bA$ is $\Delta$-modular with $\rank \bA = m$ and $\bB \subseteq \bA$, then $\overline{\bB}$ is $\Delta$-modular.}
\end{equation}
%

\section{A proof of Proposition~\ref{prop:lowerbounds}.}\label{seclowerbounds}

%
We use a generalization of the lower bound construction given by Heller.
%
Let $m,\Delta\in\mathbb{Z}_{\geq 1}$, and let $\bA \in\mathbb{Z}^{m\times n}$ consist of the following columns:

\smallskip
\begin{enumerate}[label= {(\it\roman{enumi})}]
\item \label{lb1} $\BE^i$ for every $i = 1, \dotsc, m$.
\item \label{lb2} $k \BE^1$ for every $k = 2, \dotsc, \Delta$.
\item \label{lb3} $k \BE^1 - \BE^i$ for every $k = 1, \dotsc, \Delta$ and $i = 2, \dotsc, m$.
\item \label{lb4} $\BE^i - \BE^j$ for every $2 \le i < j \le m$.
\end{enumerate}
\smallskip

The following example illustrates $\bA$ for $m=4$ and $\Delta=3$: 
\[\left[
\begin{array}{c@{\hskip .1 cm}|@{\hskip .1 cm}c@{\hskip .1 cm}|@{\hskip .1 cm}c@{\hskip .1 cm}|@{\hskip .1 cm}c@{\hskip .1 cm}!{\vrule width 1.15 pt}@{\hskip .1 cm}r@{\hskip .1 cm}|@{\hskip .1 cm}r@{\hskip .1 cm}!{\vrule width 1.15 pt}r|r|r|r|r|r|r|r|r@{\hskip .1 cm}!{\vrule width 1.15 pt}r|r|r}
1 &0 & 0 & 0 &2 & 3 &  1 & 1  & 1  &  2 & 2  & 2 & 3  & 3 & 3 & 0 & 0 & 0  \\
 0 & 1 & 0 & 0 & 0 & 0 &-1& 0  & 0  & -1 & 0  & 0 & -1 & 0 & 0 & 1 & 1 & 0  \\
 0 & 0 & 1 & 0  &0 & 0 &0  & -1 & 0  & 0  & -1 & 0 & 0  & -1& 0 & -1& 0 & 1  \\
 0 & 0 & 0 & 1  &0 & 0 &0  & 0  & -1 & 0  & 0  & -1& 0  & 0 & -1& 0 & -1&-1 
\end{array}\right].
\]
From the definition, we see that $\bA$ has differing columns, $\rank \bA=m$, and
\[
|\bA| = m+ (\Delta -1)+  (m-1)\Delta+\binom{m-1}{2} =\frac{1}{2}  \left(m^2 +m\right)+m(\Delta -1).
\]

If $m = 1$, then it is easy to verify that $\bA$ is $\Delta$-modular.
Assume that $m \ge 2$ and the proposition is true for $m-1$.
Consider a matrix $\bB \subseteq \bA$.
We prove $\Delta(\bA) \le \Delta$ by proving $|\det \bB| \le \Delta$.

Let $\hat{\bA}$ be the matrix formed by the last $m-1$ rows of $\bA $.
The matrix $\bC \subseteq \hat{\bA}$ corresponding to~\ref{lb4} form the incidence matrix of a directed graph on $m-1$ vertices, and $\hat{\bA} \setminus \bC$ is a multiset of standard unit vectors or negatives thereof. 
It is well known that $\Delta(\bC) = 1$ and $\Delta(\hat{\bA}) =1$; see~\cite[(4) on Page 268]{AS1986}.
Therefore, if $\bB$ contains a column of the form~\ref{lb2}, then by expanding $\det \bB$ along this column and using $\Delta(\hat{\bA}) =1$, we conclude $|\det \bB| \le \Delta$.

For any column of $\bA$, if we project out one of the last $m-1$ components, then the resulting column is either $\bzero$ or it is of the form of one of \ref{lb1}--\ref{lb4}, albeit in dimension $m-1$, and possibly negated.
Thus, if any of the last $m-1$ rows of $\bB$ contains exactly one non-zero entry, which necessarily equals $\pm 1$, then we expand $\det \bB$ along this row and induct on $m$ to conclude $|\det \bB| \le \Delta$.
In particular, if $\bB$ contains a column of the form~\ref{lb1}, then $|\det \bB| \le \Delta$.

Assume that $\bB$ only contains columns of the form~\ref{lb3} and~\ref{lb4}, and each of the last $m-1$ rows of $\bB$ contain at least two non-zero entries.
The invertible matrix $\bB$ must contain at least one column of the form~\ref{lb3} otherwise the first row would be all-zero.
Consider a column of the form~\ref{lb3} and suppose $\bB$ also contains a column of the form~\ref{lb4} with overlapping support, say $\bB$ contains $\ba = k \BE^1 - \BE^i$ and $\ba' = \BE^i - \BE^j$.
The matrix $[\bB |\ba + \ba'] \setminus \{\ba'\}  = [\bB |k\BE^1-\BE^j] \setminus \{\ba'\} $ has the same absolute determinant as $\bB$ and contains one more column satisfying~\ref{lb3} than $\bB$ does. 
After performing this replacement at most $m-2$ more times, we can assume that $\bB$ does not contain columns of the form~\ref{lb1} or~\ref{lb2}, each of the last $m-1$ rows of $\bB$ contains at least two non-zero entries, and $\bB$ does not contain a column of the form~\ref{lb3} and a column of the form~\ref{lb4} with overlapping supports.
Given that $\bB$ contains a column $\ba = k \BE^1 - \BE^i$ of the form~\ref{lb3} and the $i$th row of $\bB$ contains at least two non-zero entries, there exists another $\ba' \in \bB$ whose support contains the index $i$. 
After the previous replacement steps, we know that $\ba'$ must also be of the form~\ref{lb3}, that is, $\ba ' =k'\BE^1-\BE^i$ for some $k' \neq k$.
Note that $\ba - \ba' = (k-k') \BE^1$ and $[\bB |\ba - \ba'] \setminus \{\ba'\}$ contains the column $(k-k') \BE^1$, which is of the form~\ref{lb1} or~\ref{lb2}. 
Hence, $|\det \bB| = |\det([\bB |\ba - \ba'] \setminus \{\ba'\})|\le \Delta$.
\hfill\halmos

\section{Structural properties of bimodular matrices.}\label{secCircuitStructure}

In order to motivate the results in this section, we turn to a result of Heller.
Consider a TU matrix with differing columns of the form
\begin{equation}\label{eqDefnA}
\bA = 
\left[
\begin{array}{c@{\hskip .1 cm}|@{\hskip .1 cm}c}
1 & \bbeta^\intercal \\
\bzero & \hat{\bA}
\end{array}
\right],
\end{equation}
where $\bbeta \in \Z^{n-1}$ and $\hat{\bA} \in \Z^{(m-1)\times (n-1)}$.
Although $\bA$ has differing columns, the matrix $\hat{\bA}$ may not. 
After possibly multiplying columns of $\bA$ by $-1$, suppose two non-differing columns of $\hat{\bA}$ are actually equal. 
Heller showed that the set of columns in $\hat{\bA}$ with multiplicity at least two is linearly independent; see~\cite[(ii) on page 1358]{H1957}.
This linear independence is crucial in his determination of $\cc(1,m)$.
The results in this section can be viewed as a generalization of Heller's result to bimodular matrices.
It is not hard to find examples where this linear independence fails to hold for bimodular matrices. 
Rather than linear independence, we show that the set of columns in $\hat{\bA}$ with multiplicity at least two can have at most one circuit after appropriate elementary operations; see Lemma~\ref{lemStructuralCircuits}. 

This section is outlined as follows.
First, we formally define the set $M$ of columns with multiplicity at least two; see Equation~\eqref{eqM}.
Next, we provide general results of bimodular matrices in Subsection~\ref{secProperties}. 
We argue that a bimodular matrix can have at most one non-primitive column (Lemma~\ref{lem:BasicProperties}~\emph{\ref{propBasicProperties1}}), and we analyze circuits in the absence of non-primitive columns in Subsection~\ref{secPropertiesCircuits}.
Finally, we provide more precise structural statements about $\bA$ when it contains two or three linearly independent columns whose sum is divisible by two; see Subsections~\ref{subsecExtraProperties2} and~\ref{subsecExtraProperties3}.

Let $\bA \subseteq \Z^m$ be a maximal bimodular matrix with $\rank \bA = m$ and differing columns.
For any primitive column $\ba^0 \in \bA$, we can transform $\ba^0$ to $\BE^1$ via elementary operations to relabel $\bA$ as~\eqref{eqDefnA}. 
The columns in~\eqref{eqDefnA} depend on the primitive column $\ba^0$ mapped to $\BE^1$, and we make specific choices of $\ba^0$ in later subsections.
By multiplying columns of $\bA$ by $-1$, we assume that
\begin{equation}\label{eqAssumeNoNegatives}
\text{if two columns $\bb, \bc \in \hat{\bA}$ do not differ, then $\bb = \bc$.}
\end{equation}
Assumption~\eqref{eqAssumeNoNegatives} implies that $\hat{\bA}$ contains a unique maximal set of differing columns, which we denote by $\bA / \BE^1$. 
We note that if $A$ is a representation of a matroid $\mathcal{M}$, then $\bA / \BE^1$ is a representation of the matroid obtained from $\mathcal{M}$ by first contracting the element $\BE^1$, then removing ``loops'', and finally removing ``parallel'' columns that are negations of each other.
Note that the matrix $\bA / \BE^1$ may contain a column and a dilation $\alpha \bb$ for $|\alpha| \ge 2$ because we consider differing columns; this distinguishes our use of ``/\thinspace'' from the regular ``simplification" of $\bA / \BE^1$ in matroid theory.
Nevertheless, our use of the notation ``/\thinspace'' is meant to evoke the common notation for matroid contraction. 
A matrix $\bB =[\BE^1|\ba^{1}|\cdots|\ba^{m-1}] = [(1, \bzero)|(
\beta_1,\bb^{1})|\cdots|(\beta_{m-1},\bb^{m-1})] \subseteq \bA$ is a basis if and only if $\bB / \BE^1 = [\bb^{1}|\cdots|\bb^{m-1}] \subseteq \bA / \BE^1$ is a basis because $|\det \bB| = |\det \bB / \BE^1|$.
Therefore, $\bA / \BE^1$ is bimodular.
For each $\bb \in\Z^{m-1}$, we define the {\bf original set} of columns in $\bA$ corresponding to $\bb$ to be
\[
\bO(\bb) :=  \left\{
\left[\begin{array}{c}
\beta \\\bb
\end{array}
\right] \in \bA \right\}.
\]
A column $\ba \in \bA$ is said to be an {\bf original of $\bb \in \bA / \BE^1$} if $\ba = (\beta, \bb)$ for some $\beta \in \Z$.
Denote the set of columns of $\bA / \BE^1$ with  multiple originals by 
\begin{equation}\label{eqM}
\bM := \{\bb \in \bA / \BE^1 \colon |\bO(\bb)| \ge 2\}.
\end{equation}

\noindent As a reminder, throughout this section we assume bimodularity and maximality of $\bA$, as well as \eqref{eqDefnA} and~\eqref{eqAssumeNoNegatives}.

\subsection{General properties of \texorpdfstring{$\bA$}~ and \texorpdfstring{$\bM$}~.}\label{secProperties}

%
\begin{lemma}\label{lem:BasicProperties}
The matrix $\bA$ satisfies the following properties:

\smallskip
\begin{enumerate}[label= (\roman*)]
\item\label{propBasicProperties1} $\bA$ contains at most one non-primitive column, which needs to be of the form $2\ba$ for some $\ba \in \Z^m$.
Moreover, if $\bA$ only contains primitive columns, then $|\bO(\bzero)| = 1$.

\item\label{propBasicProperties3} If $\ba \in \Z^m \cap \conv [\bzero|\bA|-\bA]$, then $[\ba| \bA]$ is bimodular. 
In particular, if $\bb \in \bA/ \BE^1$ and $k := |\bO(\bb)|$, then $\bO(\bb) = \{(\beta, \bb), \dotsc, (\beta+k-1, \bb)\} \subseteq \bA$.

\item\label{propBasicProperties4} %
For each $\bb \in \bA / \BE^1$, it follows that $|\bO(\bb)| \le 3$.
\end{enumerate}
\end{lemma}
\proof{Proof.}~

\smallskip
\begin{enumerate}[label= \emph{(\roman*)}]
\item
Let $\alpha \ba \in \bA$ be non-primitive with $\alpha \ge 2$ and $\ba \in \Z$.
Let $[\alpha \ba | \bB]\subseteq \bA$ be a basis. 
We have $2 \ge |\det[\alpha \ba | \bB]| = \alpha |\det[\ba | \bB]| \ge \alpha$ because $\bA$ is bimodular and $[\ba | \bB]$ is integer valued.

%
If $2\ba, 2\bc \in \bA$ are distinct non-primitive columns, then they must be linearly independent because $\bA$ has differing columns.
Let $[2 \ba |2\bc| \bB]\subseteq \bA$ be a basis. 
We have $2 \ge |\det[2 \ba |2\bc| \bB]| \ge 4 $, which is a contradiction.
\item We can write $\ba = \sum_{i=1}^t \bc^i v_i$, where $\bc^1, \dotsc, \bc^t \in [\bzero|\bA|-\bA]$, $v_1, \dotsc, v_t \ge 0$, and $\sum_{i=1}^t v_i = 1$. 
Fix $[\ba^1|\cdots|\ba^{m-1}] \subseteq \bA$. 
By multi-linearity of the determinant and the fact that $[\bzero|\bA|-\bA]$ is bimodular, it follows that
\[
\left|\det 
\left[
\begin{array}{c@{\hskip.1cm}|@{\hskip.1cm}c@{\hskip.1cm}|@{\hskip.1cm}c@{\hskip.1cm}|@{\hskip.1cm}c}
\ba&\ba^1&\cdots&\ba^{m-1}
\end{array}
\right]\right| 
=  
\left|\sum_{i=1}^t v_i  \det\left[
\begin{array}{c@{\hskip.1cm}|@{\hskip.1cm}c@{\hskip.1cm}|@{\hskip.1cm}c@{\hskip.1cm}|@{\hskip.1cm}c}
\bc^i&\ba^1&\cdots&\ba^{m-1}
\end{array}\right] \right| \le\left|\sum_{i=1}^t v_i \right|  2= 2. 
\]
Hence, $[\ba| \bA]$ is bimodular.

Let $\bb \in \bA/ \BE^1$ with $k= |\bO(\bb)|$.
We have $\bO(\bb) = \{(\beta_1, \bb), \dotsc, (\beta_k, \bb)\} \subseteq \bA$ with $\beta_1 < \beta_2 < \cdots < \beta_k$ by Assumption~\eqref{eqAssumeNoNegatives}.
By the maximality of $A$ and the previous paragraph, we have $(\hat{\beta}, \bb) \in \bA$ for every $\beta_1 \le \hat{\beta} \le \beta_k$. 
Hence, $\{\beta_1, \beta_2, \dotsc, \beta_k\} = \{\beta_1, \beta_1+2, \dotsc, \beta_1+k-1\}$.

\item Assume $k := |\bO(\bb)| \ge 4$ for some $\bb \in \bA/ \BE^1$.
By~\emph{\ref{propBasicProperties3}}, we know that $(\beta, \bb), (\beta+3,\bb)\in \bA$ for some $\beta \in \Z$.
Let $[\bB|(\beta, \bb)| (\beta+3,\bb)] \subseteq \bA$ be a basis.
We have
\[
2 \ge 
\left|\det
\left[
\begin{array}{c|c|c}
\bB&
\begin{array}{c}
\beta \\ \bb
\end{array}&
\begin{array}{c}
\beta +3 \\ \bb
\end{array}
\end{array}
\right]
\right|
=\left|\det
\left[
\begin{array}{c|c|c}
\bB&
\begin{array}{c}
\beta \\ \bb
\end{array}&
\begin{array}{c}
3 \\ \bzero
\end{array}
\end{array}
\right]
\right|
\ge 3,
\]
which is a contradiction.
\hfill\halmos
\end{enumerate}
\endproof

\smallskip
\noindent Recall the `bar' notation defined in~\eqref{eqProject}: We write $\bC \sim ( \overline{\bC}, \bzero) \in \Z^{m \times d}$, where $ \overline{\bC} \in \Z^{\rank \bC \times d}$.

\begin{lemma}\label{lemSmallCircuits}
Assume that $\bA$ only contains primitive columns.
Let $\bC \subseteq \bM$ be a circuit.
Then

\smallskip
\begin{enumerate}[label= (\roman*)]
\item\label{propCircuit1} $3 \le |\bC| \le 4$. 
\item\label{propCircuit2} There exists a vector $\bgamma \in \Z^{|\bC|}$ such that $\left|\det \left(\bgamma^\intercal, \overline{\bC}\right)\right| = 2$, where $\left(\bgamma^\intercal, \overline{\bC}\right) \in \Z^{|\bC|\times |\bC|}$. 
\item\label{propCircuit3} If $\overline{\bC}$ is unimodular, then $\bgamma$ in~\ref{propCircuit2} satisfies $\sfrac{1}{2}\cdot \sum_{\ba \in (\bgamma^\intercal, \bC)} \ba \in \Z^m$. 
\item\label{propCircuit4} If $\overline{\bC}$ is not unimodular, then $|\bC| = 3$.
Consequently, $\bgamma$ in~\ref{propCircuit2} satisfies $\sfrac{1}{2}\cdot (\ba+\ba')$ for two columns $\ba, \ba' \in (\bgamma^\intercal, \bC)$. 
\end{enumerate}
\end{lemma}
\proof{Proof.}
Set $t := |\bC|$ and $\bC := [\bb^1|\cdots|\bb^t]$.
We have $t \ge 3$ otherwise $\bA$ contains a non-primitive column.
By Cramer's rule, $\bb^{t} = \sum_{i=1}^{t-1} \bb^i v_i$ for some $ (v_1, \dotsc, v_{t-1})=:\bv \in \{\pm 2, \pm 1, \pm \sfrac{1}{2}\}^{t-1}$.
If $v_i \in \{\pm 2\}$ for some $i$, then swap the roles of $i$ and $t$ so that $\bv \in \{\pm 1, \pm \sfrac{1}{2}\}^{t-1}$.
By~\eqref{eqPreserve}, we can assume that $\bC =  \overline{\bC}$ to simplify the remaining proof.

For each $i = 1, \dotsc, t$, the inclusion $\bb^i \in \bM$ implies $(\beta ,\bb^i) \in \bA$ for at least two choices of $\beta \in \Z$. 
By Lemma~\ref{lem:BasicProperties}~\emph{\ref{propBasicProperties3}}-\emph{\ref{propBasicProperties4}}, we have $(\beta ,\bb^i) \in \bA$ for every $\beta \in \{\beta_i, \dotsc, \beta_i+k_i\}$, where $k_i \in \{1,2\}$.
At least one value $e_i \in \{\beta_i,\dotsc, \beta_i+k_i\}$ is even. 

Define the sets
\[
\begin{array}{rcl}
\Omega &:=& \displaystyle  \left\{\bomega = (\omega_1, \dotsc, \omega_{t-1}) \colon \omega_i \in \{\beta_i, \beta_i+1\} ~~\forall ~ i = 1, \dotsc, t-1\right\}, ~\text{and}\\[.1 cm]
 \Sigma &:= & \{ \bomega^\intercal \bv \colon\bomega \in \Omega\}.
 \end{array}
\]
Each component of $\bv$ is non-zero, so $|\Sigma| \ge t \ge 3$.
For each $\bomega \in \Omega$ and $\omega_t \in\{\beta_{t},\beta_{t}+1\}$,
\[
\left| \det 
\left[ 
\begin{array}{c@{\hskip .1 cm}|@{\hskip .05 cm}c@{\hskip .05 cm}|@{\hskip .1 cm}c@{\hskip .1 cm}|@{\hskip .1 cm}c} 
\omega_1 & \multirow{2}{*}{$\cdots$} &\omega_{t-1} & \omega_t \\
\bb^1 &  & \bb^{t-1} & \bb^t
\end{array}\right]
\right| 
=
\left| \det 
\left[ 
\begin{array}{c@{\hskip .1 cm}|@{\hskip .05 cm}c @{\hskip .05 cm}|@{\hskip .1 cm}c@{\hskip .1 cm}|@{\hskip .1 cm}c} 
\omega_1 & \multirow{2}{*}{$\cdots$} & \omega_{t-1} & \omega_t - \bomega^\intercal \bv  \\
\bb^1 & & \bb^{t-1}& \bzero
\end{array}\right]
\right| =  \left|\det \left[
\begin{array}{c@{\hskip.1cm}|@{\hskip.05cm}c@{\hskip.05cm}|@{\hskip.1cm}c}
\bb^1&\cdots&\bb^{t-1}
\end{array}
\right]\right| \cdot |\omega_t - \bomega^\intercal \bv|
\le 2.
\]

Suppose $|\det [\bb^1|\cdots|\bb^{t-1}]| = 1$, i.e., $ \overline{\bC}$ is unimodular.
By Cramer's rule, $\bv \in\{\pm 1\}^{t-1}$ and $\Sigma \subseteq \Z$.
Hence, $\omega_t- \bomega^\intercal \bv \in \{\pm 2,\pm 1,0\}$ for each $\omega_t \in \{\beta_{t}, \beta_{t}+1\}$ and $\bomega^\intercal \bv \in \Sigma$.
This implies that $4 \ge |\Sigma| \ge t \ge 3$ and that there exists at least one choice $\hat{\bomega}^\intercal\bv$ and $\hat{\omega}_t$ such that $|\hat{\omega}_t- \hat{\bomega}^\intercal \bv| = 2$.
From the previous equation, we see that Property~\emph{\ref{propCircuit2}} holds with $\bgamma := (\hat{\bomega}, \hat{\omega}_t)$.
We also have that
\[
\sum_{\ba \in (\bgamma^\intercal, \bC)} \ba = \left[\begin{array}{c}\hat\omega_t\\\bb^t\end{array}\right] - \sum_{i=1}^{t-1}v_i\left[\begin{array}{c}\hat\omega_i \\\bb^i \end{array}\right]+\sum_{i=1}^{t-1} (1+v_i)\left[\begin{array}{c}\hat\omega_i \\\bb^i\end{array}\right]=  
\left[\begin{array}{c}\hat\omega_t-\hat\bomega^\intercal \bv\\ \bzero\end{array}\right]+\sum_{i=1}^{t-1} (1+v_i)\left[\begin{array}{c}\hat\omega_i \\\bb^i\end{array}\right]
\]
has all even components because $\bv \in \{\pm 1\}^{t-1}$; this implies that~\emph{\ref{propCircuit3}} holds.

Suppose $|\det [\bb^1|\cdots|\bb^{t-1}]| = 2$, i.e., $ \overline{\bC}$ is not unimodular. 
We have $\omega_t - \bomega^\intercal \bv \in \{\pm 1, \pm \sfrac{1}{2}, 0\}$ for each $\omega_t \in \{\beta_{t}, \beta_{t}+1\}$ and $\bomega^\intercal \bv \in \Sigma$.
If $t \ge 4$, then $\sigma_{\max} - \sigma_{\min} \ge \sfrac{3}{2}$,  where $\sigma_{\min}$ and $\sigma_{\max}$ are the minimum and maximum values in $\Sigma$, respectively.
Therefore, $(\beta_{t}+1)-\sigma_{\min} > 1$ or  $\beta_{t}-\sigma_{\max} <-1$, which is a contradiction.
Thus, $|\Sigma| = t = 3$, and $\bv\in\{\pm\sfrac12\}^{t-1}$. 
Furthermore, there is at least one choice $\hat{\bomega}^\intercal \bv$ and $\hat{\omega}_t$ such that $|\hat{\omega}_t- \hat{\bomega}^\intercal \bv| = 1$. 
Property~\emph{\ref{propCircuit2}} holds with $\bgamma := (\hat{\bomega}, \hat{\omega}_t)$.
We have $\bb^3 = \bb^1v_1 +\bb^2v_2 \in \Z^{m-1}$, which implies $\sfrac{1}{2}\cdot(\bb^1+\bb^2) \in \Z^{m-1}$ because $\bv\in\{\pm \sfrac12\}^2$.
Similarly, $\sfrac{1}{2}\cdot(\hat{\omega}_1+\hat{\omega}_2) \in \Z$ because $\bv\in\{\pm \sfrac12\}^2$, $|\hat{\omega}_3- \hat{\bomega}^\intercal \bv| = 1$, and $\hat{\omega}_3 \in \Z$. 
Hence, $\sfrac{1}{2} \cdot ((\hat{\omega}_1,\bb^1)+(\hat{\omega}_2,\bb^2)) \in \Z^m$, which implies that~\emph{\ref{propCircuit4}} holds.
\hfill\halmos
\endproof

\subsection{Circuits in \texorpdfstring{$\bM$}~ when \texorpdfstring{$\bA$}~ contains only primitive columns.}\label{secPropertiesCircuits}

In this subsection, we assume that $\bA$ only contains primitive columns and $\bM$ contains a circuit; these assumptions allow us to apply Lemma~\ref{lemSmallCircuits}.
Choose $\bB^*\subseteq \bA$ satisfying
\begin{equation}\label{eqMinimalAssumption}
\text{$\bB^*$ is linearly independent and $\frac12\cdot\sum_{\ba\in \bB^*}\ba\in\Z^m$},
\end{equation}
and minimizing $|\bB^*|$.
The set $\bB^*$ exists and $|B^*| \le 4$ by Lemma~\ref{lemSmallCircuits}~\emph{\ref{propCircuit3}}-\emph{\ref{propCircuit4}}.
Furthermore, $2 \le |\bB^*|$ because we assumed that $A$ only contains primitive columns.
After applying elementary operations, we assume that
\begin{equation}\label{eqBForm}
\bB^*=\left[
\begin{array}{c@{\hskip .1 cm}|@{\hskip .1 cm}c@{\hskip .1 cm}|@{\hskip .1 cm}c@{\hskip .1 cm}|@{\hskip .1 cm}c}
\BE^1 
&
\cdots
&
\BE^{|\bB^*|-1}
&
\BE^1 +\cdots+\BE^{|\bB^*|-1}+2\BE^{|\bB^*|}
\end{array}\right].
\end{equation}
%


\begin{lemma}\label{lemSpan}
%
%
%
After possibly multiplying columns of $\bA$ by $-1$, we can assume that
\begin{equation}\label{BSpan} 
\bA\cap \Span \bB^* \subseteq \bB^* \cup \left\{\BE^{|\bB^*|}+\sum_{i=1}^{|\bB^*|-1} \alpha_i\BE^i
\colon\alpha_1,\dotsc,\alpha_{|\bB^*|-1}\in\{0,1\}\right\}.
\end{equation}
Furthermore, $|\bO(\bb)| =1$ for each $\bb \in \bB^* / \BE^1$ and 
\begin{equation}\label{M=2} 
M= \{\bb \in \bA/ \BE^1 \colon |\bO(\bb)| = 2\}. 
\end{equation}
\end{lemma}

\proof{Proof.}
Set $s := |\bB^*| \ge 2$.
Let $\ba \in \bA\cap \Span \bB^* \setminus \bB^*$.
By Cramer's rule and the assumption that $\bA$ is bimodular, we can write $\ba = \bB^* \bv$, where $\bv = (v_1, \dotsc, v_{s}) \in  \{ \pm 1, \pm \sfrac{1}{2}, 0\}^{s}$. 
For proving \eqref{BSpan}, it suffices to show $v_{s}=\sfrac{1}{2}$ because $\ba = \bB^* \bv \in\mathbb{Z}^m$ implies $v_{i}\in\{\pm\sfrac{1}{2}\}$ for all $i=1,\dotsc,s-1$. 

Set $I:=\supp (v_1,\dotsc,v_{s-1})$. 
Suppose $v_{s}=0$; then $\bv \in \{ 0,\pm 1 \} ^{s}$ and $|I| \le 1$ otherwise $\sfrac{1}{2} \cdot (\ba + (\BE^1+\dots+\BE^{s-1}+2\BE^{s})+\sum_{i\in \{ 1,\dotsc,s-1 \}\setminus I} \BE^i ) \in \Z^m$, which contradicts the minimality of $\bB^*$. 
However, $|I|\leq 1$ implies $\ba \in \bB^*$, which is a contradiction. 
Suppose $v_{s} \in\{\pm 1\}$; then $\bv \in \{ 0,\pm 1 \} ^{s}$ and $|I|=0$ otherwise $\sfrac{1}{2} \cdot (\ba + \sum_{i\in \{ 1,\dotsc,s-1\}\setminus I} \BE^i ) \in \Z^m$, which contradicts the minimality of $\bB^*$. 
However, $|I|=0$ implies $\ba = \bzero$, which contradicts that $A$ has differing columns. 
Thus, $v_{s} = \pm \sfrac{1}{2}$ and by possibly replacing $\ba$ by $-\ba$, we assume that $v_{s} = \sfrac{1}{2}$.
This proves \eqref{BSpan}.

It follows directly from~\eqref{BSpan} that $|\bO(\bb)| =1$ for each $\bb \in \bB^* / \BE^1$.
Assume to the contrary that $|\bO(\bb)| \ge 3$ for some $\bb \in \bA/ \BE^1$. 
It follows from Lemma \ref{lem:BasicProperties} \emph{\ref{propBasicProperties3}}-\emph{\ref{propBasicProperties4}} that $\bO(\bb)= \{\ba, \ba+\BE^1,\ba+2\BE^1\}$ for some $\ba \in \bA$.
Inclusion~\eqref{BSpan} implies that if $\bc \in \bA \cap \Span \bB^*$, then $\bc + 2 \BE^1 \not \in \bA$.
Therefore $\ba \not \in \Span \bB^* $, and $\bC = [\BE^2 |\cdots|\BE^{s-1} |\BE^1+\cdots+\BE^{s-1} + 2\BE^{s}|\ba|\ba+2\BE^1] \subseteq \bA$ has linearly independent columns.
Recall~\eqref{eqProject}: we have $\bC \sim ( \overline{\bC}, \bzero)$, where 
\[
 \overline{\bC} = \left[
\begin{array}{c@{\hskip .1 cm}|@{\hskip .1 cm}c@{\hskip .1 cm}|@{\hskip .1 cm}c@{\hskip .1 cm}|@{\hskip .1 cm}c@{\hskip .1 cm}|@{\hskip .1 cm}c@{\hskip .1 cm}|@{\hskip .1 cm}c}
\BE_{s+1}^2 
&
\cdots
&
\BE_{s+1}^{s-1} 
&
\BE_{s+1}^1+\cdots+\BE_{s+1}^{s-1} + 2\BE_{s+1}^{s}
&
\overline{\ba}
&
\overline{\ba}+2\BE_{s+1}^1 
\end{array}\right]\in\mathbb{Z}^{(s+1)\times (s+1)}
\]
is invertible and $\left|\det  \overline{\bC}\right|\geq 4$, which contradicts~\eqref{eqPreserve}.
\hfill \halmos
\endproof

\smallskip
Define the matrices
\begin{equation}\label{eqn:C0}
\bC^* :=\left\{
\BE^{|\bB^*|-1}_{m-1}+\sum_{i=1}^{|B^*|-2}\alpha_i\BE^i_{m-1}
:~\alpha_1,\dotsc,\alpha_{|\bB^*|-2}\in\{0,1\}\right\}
\end{equation}
and
\begin{equation}\label{eqn:D}
\bD^* := \bB^*/\BE^1 =\left[
\begin{array}{c@{\hskip .1 cm}|@{\hskip .1 cm}c@{\hskip.1 cm}|@{\hskip .1cm}c@{\hskip .1 cm}|@{\hskip.1cm}c@{\hskip .1 cm}|@{\hskip.1cm}c}
\BE^1_{m-1}&\BE^2_{m-1}&\cdots&\BE^{|\bB^*|-2}_{m-1} & \BE^1_{m-1}+\cdots+\BE^{|\bB^*|-2}_{m-1}+2\BE^{|\bB^*|-1}_{m-1}
\end{array}\right].
\end{equation}
By \eqref{BSpan}, we can assume that $\bC^*$ contains all columns in $(\bA\cap\Span \bB^*)/\BE^1$ that have multiple originals, i.e., $\bM\cap \Span \bD^* \subseteq \bC^*$.
\begin{lemma}\label{MpropertyPair}
If $[\bb|\bb+\bd]\subseteq \bM$ for some $\bb\in \bM\setminus \bC^*$ and $\bd\in \Span \bD^*$, then $\bd\in [\bD^*|-\bD^*]$.
\end{lemma}
\proof{Proof.}
Set $s := |\bB^*|-1 \ge 1$.
By Cramer's rule and the bimodularity of $\bA/ \BE^1$, we have $\bd=\bD^*\bv$, where $\bv=(v_1,\dots,v_s) \in \{\pm 1, \pm \sfrac{1}{2}, 0\}^s$.
Set $\bD^* = [\bd^1|\cdots|\bd^s]$.
Recall \eqref{eqProject}:
\[
[\bD^*|\bb|\bb+\bd]\sim 
\left[\begin{array}{c@{\hskip.1 cm}|@{\hskip.1 cm}c@{\hskip.1 cm}|@{\hskip.1 cm}c}\overline{\bD^*}&\bzero&\overline{\bd}\\0&1 &1\\
\bzero & \bzero & \bzero 
\end{array}\right]
\]
where $\overline{\bD^*}=[\overline{\bd^1}| \cdots|\overline{\bd^s}]=[\BE^1_s|\cdots|\BE^{s-1}_s|\BE^1_s+\cdots+\BE^{s-1}_s+2\BE^s_s]$ and $\overline{\bd}= \overline{\bD^*}\bv$.
For every choices of $(\gamma_1,\bb)$ and  $(\gamma_2,\bb+\bd)$ in $\bA$, we have
\[
\left|\det
\left[\begin{array}{c@{\hskip.1 cm}|@{\hskip.1 cm}c@{\hskip.1 cm}|@{\hskip.1 cm}c@{\hskip.1 cm}|@{\hskip.1 cm}c@{\hskip.1 cm}|@{\hskip.1 cm}c@{\hskip.1 cm}|@{\hskip.1 cm}c}
0 & & 0 &  1 & \gamma_1 & \gamma_2\\
\bd^1 &\cdots & \bd^{s-1} & \bd^s & \bzero & \overline{\bd}\\
0 & & 0 & 0 & 1 & 1
\end{array}
\right]\right| 
=2|\gamma_2-\gamma_1 - v_s|\le 2.
\]

If $v_s\not\in \Z$, then $2(\gamma_2 - \gamma_1-v_s)$ is odd and contained in $ \{\pm 2, \pm 1, 0\}$.
There are two choices for both $\gamma_1$ and $\gamma_2$ because $[\bb|\bb+\bd] \subseteq \bM$.
However, this means that there are at least three distinct odd values of $2(\gamma_2 - \gamma_1-v_s)$ in $ \{\pm 2, \pm 1, 0\}$, which is a contradiction.
Hence, $v_s\in \Z$.
This implies that $\bv \in\{\pm1,0\}^s$ because $\bD^* \bv = \bd \in \Z^{m-1}$.

Set $I:=\supp(v_1,\dotsc,v_{s})$.
If $|I| \ge 2$, then $\sfrac{1}{2}\cdot  (\bd + \sum_{i\in \{1,\dots,s\}\setminus I}\bd^i)\in \Z^{m-1}$.
This implies there exist originals of $\bb,\bb+\bd,$ and $\bd^i$ for each $i\in \{1,\dots,s\}\setminus I$ that satisfy \eqref{eqMinimalAssumption}.
However, there are only $s+2-|I|<s+1=|\bB^*|$ columns here, which contradicts the minimality of $\bB^*$.
Hence, $|I|=1$ and $\bd\in [\bD^*|-\bD^*]$.
\hfill \halmos
\endproof

\begin{lemma}\label{MpropertyLinInd}
If $\bC=[\bb^{1}|\cdots|\bb^{t}]\subseteq \bM$ is a circuit, then $[\bD^*|\bb^{j_1}|\cdots|\bb^{j_{t-1}}]$ contains a circuit for every choice of indices $j_1, \dotsc, j_{t-1} \in \{1, \dotsc, t\}$.
\end{lemma}
\proof{Proof.}
Set $s:= |\bB^*|-1 \ge 1$, and set $\bD^* = [\bd^1|\cdots|\bd^s]$.
Assume to the contrary that $[\bD^*|\bb^{1}|\cdots|\bb^{{t-1}}]$ has linearly independent columns.
Using linear independence and \eqref{eqProject}, we have 
\[
[\bD^*|\bC] \sim \left[\begin{array}{c@{\hskip .1 cm}|@{\hskip .1 cm}c}
\overline{\bD}& \bzero\\\bzero& \overline{\bC} \\ \bzero & \bzero\end{array}\right],
\]
where $\overline{\bD^*}=[\overline{\bd^1}|\cdots|\overline{\bd^s}] = [\BE^1_s|\cdots|\BE^{s-1}_s|\BE^1_s+\cdots+\BE^{s-1}_s+2\BE^s_s]$ and $ \overline{\bC}\in\Z^{(t-1)\times t}$.
By Lemma \ref{lemSmallCircuits} \emph{\ref{propCircuit2}}, there exists some $\bgamma = (\gamma_1, \dotsc, \gamma_t) \in \Z^t$ such that $(\bgamma^\intercal, \bC) \subseteq \bA$ and $|\det (\bgamma^\intercal,  \overline{\bC})| = 2$.
Therefore, 
\[
\left|\det
\left[
\begin{array}{c@{\hskip.1 cm}|@{\hskip.1 cm}c@{\hskip.1 cm}|@{\hskip.1 cm}c@{\hskip.1 cm}|@{\hskip.1 cm}c@{\hskip.1 cm}|@{\hskip.1 cm}c}
0 &  & 0 & 1 &\bgamma^\intercal \\
\overline{\bd^1} & \cdots & \overline{\bd^{s-1}} & \overline{\bd^s} & \bzero \\
0 &  & 0 & 0 &  \overline{\bC}
\end{array}
\right]\right| 
=\left|\det \overline{\bD^*}\right|\cdot\left|\det (\bgamma^\intercal,  \overline{\bC})\right|=4,
\]
which contradicts~\eqref{eqPreserve}.
\hfill \halmos
\endproof
\begin{lemma}\label{MpropertyCircuit}
If $\bC\subseteq \bM$ is a circuit and $|\bC\setminus \bC^*|\ge 2$, then $\bC\setminus \bC^*=[\bb|\bb+\bd]$ for some $\bb\in \bM\setminus \bC^*$ and $\bd\in \bD^*$.
Given that $|\bC| \in \{3,4\}$ from Lemma~\ref{lemSmallCircuits}~{\ref{propCircuit1}}, it follows that $\bC \cap \bC^* \neq \emptyset$.
\end{lemma}
\proof{Proof.}
Set $s:= |\bB^*|-1 \ge 1$ and $t:= |C|$.
By Lemma \ref{MpropertyLinInd}, we know that $\rank[\bD^*|\bC] \le \rank \bD^*+\rank \bC-1=s+t-2$.
Also, by $|\bC\setminus \bC^*|\ge 2$ and $\bM\cap\Span \bD^*\subseteq \bC^*$, we know that $\rank[\bD^*|\bC]\ge \rank \bD^* +1 = s+1$.
By Lemma \ref{lemSmallCircuits}~\emph{\ref{propCircuit1}}, we have $t\in\{3,4\}$.
In both cases, we argue that $|C \setminus C^*| = 2$ and $\rank[\bD^*|\bC] = s+1$.
%
It will then follow from Lemma~\ref{MpropertyPair} and after possibly multiplying the column by $-1$, that $\bC\setminus \bC^*=[\bb|\bb+\bd]$ for $\bb\in \bM\setminus \bC^*$ and $\bd\in \bD^*$.

Assume that $t=3$; then $\rank[\bD^*|\bC] = s+1$. 
If $|\bC\setminus \bC^*|=3$, then $[\bD^*|\bC]\sim [\bD^*|\BE^{s+1}_{m-1}|\BE^{s+1}_{m-1}+\bd^1|\BE^{s+1}_{m-1}+\bd^2]$ for distinct $\bd^1,\bd^2\in\Span \bD^*$.
By Lemma~\ref{MpropertyPair}, we have $\bd^1,\bd^2\in [\bD^*|-\bD^*]$.
The matrix $[\BE^{s+1}_{m-1}+\bd^1|\BE^{s+1}_{m-1}+\bd^2]$ is contained in $\bM$ but $\bd^2-\bd^1\notin [\bD^*|-\bD^*]$ for any two distinct columns in $D^*$; this contradicts Lemma~\ref{MpropertyPair}.
Therefore, $|\bC\setminus \bC^*|=2$ when $t = 3$.

Assume that $t=4$.
If $\rank[\bD^*|\bC] = s+1$, then $|\bC\setminus \bC^*|=2$ as in the case $t=3$.
Assume to the contrary that $\rank[\bD^*|\bC] = s+2$.
By \eqref{eqProject}, we can assume that
\[
[\bD^*|\bC]\sim
\left[\begin{array}{c@{\hskip .1 cm}|@{\hskip .1 cm}c@{\hskip .1 cm}|@{\hskip .1 cm}c@{\hskip .1 cm}|@{\hskip .1 cm}c}\overline{\bD^*}& \bzero &\overline{\bD^*}\bu^3 &\overline{\bD^*}\bu^4\\\bzero& \mathbb{I}_2&\bv^3&\bv^4\\
\bzero & \bzero & \bzero & \bzero \end{array}\right],
\]
where $\overline{\bD^*} =[\BE^1_s|\cdots|\BE^{s-1}_s|\BE^1_s+\cdots+\BE^{s-1}_s+2\BE^s_s]$, $\bu^3 = (u^3_1, \dotsc, u^3_s)$ and $\bu^4$ are contained in $\{0,\pm\sfrac12,\pm1\}^s$, and $\bv^3 = (v^3_1, v^3_2)$ and $\bv^4 = (v^4_1, v^4_2)$ are contained in $\{0,\pm1\}^2$.
Lemma~\ref{lemSmallCircuits}~\emph{\ref{propCircuit4}} implies $\bC$ is unimodular because $t = 4$, and Lemma~\ref{lemSmallCircuits}~\emph{\ref{propCircuit3}} implies $\sfrac{1}{2} \cdot (\bone +\bv^3+\bv^4) \in \Z^{2}$.
We derive a contradiction in two cases.

First, assume that $\bv^3\in\{\pm1\}^2$ or $\bv^4\in \{\pm1\}^2$.
Say $\bv^3\in\{\pm1\}^2$; then $\bv^4=\bzero$ and $|\bC\setminus \bC^*|=3$.
For $\delta \in \Z$ and $\bgamma = (\gamma_1, \gamma_2) \in 
\Z^2$, the matrix
\[
\bE(\bgamma, \delta) := 
\left[
\begin{array}{c@{\hskip .1 cm}|@{\hskip .1 cm}c@{\hskip .1 cm}|@{\hskip .1 cm}c@{\hskip .1 cm}|@{\hskip .1 cm}c@{\hskip .1 cm}|@{\hskip .1 cm}c@{\hskip .1 cm}|@{\hskip .1 cm}c}
0 &  & 0 & 1 & \bgamma^\intercal & \delta\\ 
\overline{\bd^1} & \cdots & \overline{\bd^{s-1}} & \overline{\bd^s} & \bzero &\overline{\bD^*}\bu^3\\
\bzero& & \bzero&\bzero&\mathbb{I}_2& \bv^3\end{array}\right] \in \Z^{(s+3)\times (s+3)}
\]
has an absolute determinant of $2|\delta - \bgamma^\intercal \bv^3 - u^3_s| = |2(\delta- \bgamma^\intercal \bv^3) - 2u^3_s| \in \{0,1,2\}$.
Given that $\bC \subseteq \bM$, there are two choices for each of $\gamma_1, \gamma_2,$ and $\delta$ such that $(\bE(\bgamma, \delta), \bzero) \subseteq \bA$.
Thus, there are at least four distinct values of $2(\delta- \bgamma^\intercal \bv^3) - 2u^3_s$ in $\{\pm 2, \pm 1, 0\}$ that have the same parity, namely, the same parity as $2u^3_s$.
However, this is a contradiction. 

Second, assume that $\bv^3\in\{\pm\BE^1_2\}$ and $\bv^4\in\{\pm\BE^2_2\}$ or that $\bv^3\in\{\pm\BE^2_2\}$ and $\bv^4\in\{\pm\BE^1_2\}$; then $|\bC\setminus \bC^*|=4$.
By Lemma \ref{MpropertyPair}, we have $\overline{\bD^*}\bu^3, \overline{\bD^*}\bu^4 \in [\overline{\bD^*}|-\overline{\bD^*}]$.
By possibly multiplying the column by $-1$, we assume $\overline{\bD^*}\bu^3\in\overline{\bD^*}$.
Set $\bF := \bD^* \setminus \{\bD^* \bu^3\}$ and $\overline{\bF} := \overline{\bD^*} \setminus \{\overline{\bD^*}\bu^3\}$.
By Lemma \ref{lemSmallCircuits} \emph{\ref{propCircuit2}}, there exists a vector $\bgamma =(\gamma_1, \dotsc, \gamma_4) \in \Z^4$ such that $(\bgamma^\intercal, \bC) \subseteq \bA$ and $2 = \left|\det (\bgamma^\intercal,  \overline{\bC})\right| $.
Let $\bdelta \in \Z^{s-1}$ be such that $(\bdelta^\intercal, \bF) \subseteq \bA$.
We have
\[
\left|
\det
\left[
\begin{array}{c@{\hskip .1 cm}|@{\hskip .1 cm}c@{\hskip .1 cm}|@{\hskip .1 cm}c@{\hskip .1 cm}|@{\hskip .1 cm}c@{\hskip .1 cm}|@{\hskip .1 cm}c}
\bdelta^\intercal &\gamma_1 & \gamma_2 & \gamma_3& \gamma_4\\
\overline{\bF} & \bzero & \bzero & \overline{\bD^*}\bu^3 & \overline{\bD^*}\bu^4 \\
\bzero& \BE^1_2&\BE^2_2 & \bv^3&\bv^4 
\end{array}\right]
\right|
=
\left|
\det
\left[
\begin{array}{c@{\hskip .1 cm}|@{\hskip .1 cm}c@{\hskip .1 cm}|@{\hskip .1 cm}c}
\overline{\bF} & \bzero  & \overline{\bD^*}\bu^3\\
\bzero& \mathbb{I}_2 & \bv^3
\end{array}\right]
\right|\cdot\left|\det (\bgamma^\intercal,  \overline{\bC})\right| = 2 \left|\det\overline{\bD^*}\right|= 4,
\]
which is a contradiction.

Therefore, $|\bC\setminus \bC^*|=2$ and $\rank[\bD^*|\bC] = s+1$ when $t = 4$.
\hfill \halmos
\endproof

\smallskip
We arrive at our main result in this section.
We repeat assumptions for the reader. 

\begin{lemma}\label{lemStructuralCircuits}
Assume that $\bA$ only contains primitive columns and let $\bC^1 \subseteq \bM$ be a circuit. 
Choose $\bB^*\subseteq \bA$ that satisfies~\eqref{eqMinimalAssumption} and minimizes $|\bB^*|$, and assume that $\bB^*$ has the form~\eqref{eqBForm}.
Recall $\bC^*$ from~\eqref{eqn:C0} and $\bD^*$ from~\eqref{eqn:D}.
\smallskip%
\begin{enumerate}[label= (\roman*)]
\item\label{struc1} If $|\bB^*| = 2$, then $\bC^1 = [\bC^*|\bb|\bb+2\BE^1_{m-1}]$ for some $\bb \in \bA/ \BE^1$.
Furthermore, $\bC^1$ is the unique circuit in $\bM$, so $|\bM| \le m$.

\item\label{struc2} If $|\bB^*| = 3$, then $\bC^1 = [\bC^*|\bb|\bb+\bd]$ for some $\bb\in \bA/\BE^1$ and $\bd\in \bD^*$.
Furthermore, $\bC^1$ is the unique circuit in $\bM$, so $|\bM| \le m$.

\item\label{struc3} If $|\bB^*| = 4$, then either $\bC^1 = \bC^*$, or $|\bC^1 \cap \bC^*| = 2$ and $|\bC^1| = 4$.
Moreover, if $\bM$ contains multiple circuits, say $\bC^1 \neq \bC^*$, then $\bM$ contains precisely three circuits: $\bC^1$, $\bC^*$, and the symmetric difference $\bC^1 \triangle \bC^*$. 
Regardless of the number of circuits, $|\bM|\le m+1$.
\end{enumerate}
\end{lemma}

\proof{Proof.}
We have $|\bC^1| \in \{3,4\}$ by Lemma~\ref{lemSmallCircuits}~\emph{\ref{propCircuit1}}.

\smallskip
\begin{enumerate}[label= (\emph{\roman*})]
\item 
Note that $\bC^* = [\BE^1_{m-1}]$, so $|\bC^1\setminus \bC^*|\ge 2$.
It follows from Lemma~\ref{MpropertyCircuit} that $|\bC^1\setminus \bC^*|=2$ and $\bC\setminus \bC^* = [\bb|\bb+2\BE^1_{m-1}]$ for some $\bb\in \bM\setminus \bC^*$.
Therefore, $\bC^1=[\bC^*|\bb|\bb+2\BE^1_{m-1}]$.

If $\bM$ contains another circuit $\bC^2$, then $\bC^2 = [\bC^*|\bb'|\bb'+2\BE^1_{m-1}]$ with $\bb \neq \bb'$.
The column $\BE^1_{m-1}$ is linearly dependent on $[\bb|\bb+2\BE^1_{m-1}]$ and on $[\bb'|\bb'+2\BE^1_{m-1}]$.
Hence, $[\bb|\bb+2\BE^1_{m-1}|\bb'|\bb'+2\BE^1_{m-1}] \subseteq \bM\setminus \bC^*$ contains a circuit $\bC^3$.
However, $|\bC^3 \setminus \bC^*| = |\bC^3| \ge 3$ because $\bC^3 \subseteq \bM\setminus \bC^*$; this contradicts Lemma~\ref{MpropertyCircuit}.

\item Note that $\bC^* = [\BE^2_{m-1}|\BE^1_{m-1}+\BE^2_{m-1}]$, so $|\bC^1\setminus \bC^*|\ge 1$.
If $|\bC^1\setminus \bC^*| = 1$, then $\bC^1\subseteq \Span \bC^*=\Span \bD^*$, which contradicts $\bM\cap\Span \bD^* \subseteq \bC^*$.
Thus, $|\bC^1\setminus \bC^*|\ge 2$.
It follows from Lemma~\ref{MpropertyCircuit} that $|\bC^1\setminus \bC^*|=2$ and $\bC^1\setminus \bC^* = [\bb|\bb+\bd]$ for some $\bb\in \bM\setminus \bC^*$ and $\bd\in D^*$.
Hence, $\bC^1=[\bC^*|\bb|\bb+\bd]$ because $\rank[\BE^2_{m-1}|\bb|\bb+\bd]=\rank[\BE^1_{m-1}+\BE^2_{m-1}|\bb|\bb+\bd]=3$.

If $\bM$ contains another circuit $\bC^2$, then $\bC^2 = [\bC^*|\bb'|\bb'+\bd']$ for some $\bb'\in \bM \setminus\bC^*$ and $\bd'\in \bD^*$.
We know that $[\bb|\bb+\bd|\bb'|\bb'+\bd']\subseteq \bM\setminus \bC^*$ does not contain a circuit because such a circuit would not satisfy Lemma~\ref{MpropertyCircuit}.
Therefore, $\rank[\bb|\bb+\bd|\bb'|\bb'+\bd'] =4$. 
However, $\BE^1_{m-1}+\BE^2_{m-1}$ is linearly dependent on  $[\BE^2_{m-1}|\bb|\bb+\bd]$ and on $[\BE^2_{m-1}|\bb'|\bb'+\bd']$, which implies that $[\BE^2_{m-1}|\bb|\bb+\bd|\bb'|\bb'+\bd']$ is a circuit with five columns; this contradicts Lemma~\ref{lemSmallCircuits}~\emph{\ref{propCircuit1}}.

\item If $|\bC^1| = 3$ and $\bC^1$ is unimodular, then $|\bB^*| \le 3$ by Lemma~\ref{lemSmallCircuits}~\emph{\ref{propCircuit3}}, which contradicts the minimality of $|\bB^*|$ in Case~\emph{\ref{struc3}}.
If $|\bC^1| = 3$ and $\bC^1$ is not unimodular, then $|\bB^*| \le 2$ by Lemma~\ref{lemSmallCircuits}~\emph{\ref{propCircuit4}}, which again contradicts the minimality of $|\bB^*|$ in this case.
Hence, $|\bC^1| = 4$.

Suppose $\bC^1 \neq \bC^*$; then $|\bC^1\setminus \bC^*|\ge 1$ because both matrices are circuits.
%
%
Recall $M \cap \Span D^* \subseteq C^*$; see the sentence after~\eqref{eqn:D}.
If $|\bC^1\setminus \bC^*| = 1$, then $\bC^1\subseteq M \cap \Span \bC^*=M \cap \Span D^* \subseteq C^*$ because $\bC^1$ is a circuit; this contradicts $|\bC^1\setminus \bC^*| = 1$.
It then follows from Lemma~\ref{MpropertyCircuit} that $|\bC^1\setminus \bC^*|=2$.
Hence, $\bC^1 = [\bc^1|\bc^2|\bb^1|\bb^1+\bd^1]$ for some $\bc^1, \bc^2 \in \bC^*$, $\bb^1\in \bM\setminus \bC^*$, and $\bd^1\in \bD^*$.
We have
\begin{equation}\label{eqBigMatrix1}
[\bC^*|\bb^1|\bb^1+\bd^1]
\sim 
\left[ 
\begin{array}{c@{\hskip .1 cm}|@{\hskip .1 cm}c@{\hskip .1 cm}|@{\hskip .1 cm}c@{\hskip .1 cm}|@{\hskip .1 cm}c@{\hskip .1 cm}|@{\hskip .1 cm}c@{\hskip .1 cm}|@{\hskip .1 cm}c}
0 & 1 & 0 & 1 & 0 & d^1_1\\
0 & 0 & 1 & 1 & 0 & d^1_2\\
1 & 1 & 1 & 1 & 0 & d^1_3\\
0 & 0 & 0 & 0 & 1 & 1\\
\bzero & \bzero & \bzero & \bzero& \bzero& \bzero
\end{array}
\right],
\end{equation}
where $\overline{\bd} := (d^1_1, d^1_2, d^1_3) \in \overline{\bD^*} = [(1,0,0)| (0,1,0)|(1,1,2)]$.

Given that $\bC^1$ is a circuit, $\bd^1$ linearly depends on $\bc^1$ and $\bc^2$.
From this and~\eqref{eqBigMatrix1}, we can determine $[\bc^1|\bc^2]$ (the index sets refer to the matrix on the right hand side of~\eqref{eqBigMatrix1}):
\begin{equation}\label{eqDeterminec1c2}
\begin{array}{l}
\text{If $\overline{\bd} = (1,0,0)$, then $[\bc^1| \bc^2]$ is indexed by $\{1, 2\}$ or $\{3,4\}$.}\\
\text{If $\overline{\bd} = (0,1,0)$, then $[\bc^1| \bc^2]$ is indexed by $\{1,3\}$ or $\{2,4\}$.}\\
\text{If $\overline{\bd} = (1,1,2)$, then $[\bc^1| \bc^2]$ is indexed by $\{1, 4\}$ or $\{2,3\}$.}
\end{array}
\end{equation}

\begin{claim}\label{claime1e2e3e4_3}
If $\bC^2 \subseteq \bM$ is a circuit and $\bC^2 \neq \bC^*$, then $\bC^1\setminus \bC^* = \bC^2 \setminus \bC^*$.
\end{claim}
\proof{Proof.}
Assume to the contrary that there exists a circuit $\bC^2 \subseteq \bM$ with $\bC^2 \neq \bC^*$ and $\bC^1\setminus \bC^* \neq \bC^2 \setminus \bC^*$.
As was the case with $\bC^1$, we can write $\bC^2= [\bc^3|\bc^4|\bb^2|\bb^2+\bd^2]$ for some $\bc^3, \bc^4 \in \bC^*, \bb^2\in \bM\setminus \bC^*$ and $\bd^2\in \bD^*$.
We know that $[\bb^1|\bb^1+\bd^1]\ne [\bb^2|\bb^2+\bd^2]$ because $\bC^1\setminus \bC^* \ne \bC^2 \setminus \bC^*$.
Also, $[\bb^1|\bb^1+\bd^1|\bb^2|\bb^2+\bd^2]\subseteq \bM\setminus \bC^*$ does not contain a circuit because such a circuit would not satisfy Lemma~\ref{MpropertyCircuit}.
Thus, $\rank [\bb^1|\bb^1+\bd^1|\bb^2|\bb^2+\bd^2] = 4$ and $\bd^1 \neq \bd^2$.
It follows from~\eqref{eqDeterminec1c2} that $\{\bc^1,\bc^2\} \setminus \{\bc^3,\bc^4\}\neq \emptyset $; say $\bc^1 \not \in \{\bc^3,\bc^4\}$. 
Similarly, we can assume that $\bc^3 \not \in \{\bc^1,\bc^2\}$.
Consequently, any three columns of $[\bd^1|\bd^2|\bc^1|\bc^3]$ are linearly independent.
However, this implies $[\bb^1|\bb^1+\bd^1|\bb^2|\bb^2+\bd^2|\bc^1|\bc^3] $ is a circuit of cardinality six, which contradicts Lemma \ref{lemSmallCircuits} \emph{\ref{propCircuit1}}.
\hfill$\diamond$
\endproof

\smallskip

Recall $|\bM\cap \bC^*|\ge 2$.
By Claim~\ref{claime1e2e3e4_3}, every circuit in $\bM$ is contained in~\eqref{eqBigMatrix1}.
By Lemma \ref{MpropertyCircuit}, any circuit in $\bM$ besides $\bC^*$ uses columns 5 and 6 alongside two of the first four columns. 
This shows that either $\bM$ only contains one circuit, namely $\bC^1$, or $\bM$ contains the three circuits $\bC^1$, $\bC^*$, and $\bC^1 \triangle \bC^*$.

Suppose $|\bM\cap \bC^*|\le 3$.
Hence, $\bC^*$ is not contained in $\bM$.
For any two columns in $\bM\cap \bC^*$ and any choice of $(d^1_1, d^1_2, d^1_3)$, there exists at most one pair of columns in $\bM\cap \bC^*$ that form a circuit with columns 5 and 6; see~\eqref{eqDeterminec1c2}.
Hence, $\bM$ contains at most one circuit, so $|\bM| \le  m$.

Suppose $|\bM\cap \bC^*|=4$; then $\bM \cap \bC^* = \bC^*$.
Let $\bc \in \bM\cap \bC^*$.
The matrix $\bM\setminus \{\bc\}$ satisfies $|\bM \cap \bC^* \setminus \{\bc\}| = 3$.
Therefore, $|\bM\setminus \{\bc\}| \le m$ from the previous paragraph. 
Hence, $|\bM| \le m+1$.
\hfill\halmos
\endproof
\end{enumerate}
%

\subsection{Additional structural properties when \texorpdfstring{$|\bB^*| = 2$}~.}\label{subsecExtraProperties2}
%

\begin{lemma}\label{lem4ptlineplus3}
Assume that $\bA$ only contains primitive columns. 
Choose $\bB^*\subseteq \bA$ that satisfies~\eqref{eqMinimalAssumption} and minimizes $|\bB^*|$, and assume that $\bB^*$ has the form~\eqref{eqBForm}.
If $|\bB^*| = 2$, $m\ge 3$ and $|\bM|=m$, then
\begin{equation}\label{eqn4ptlineplus3}
 |\bA|\le \frac{1}{2} \left(m^2+m\right) +3. 
\end{equation}
\end{lemma}

\proof{Proof.}
%

%
%


The assumption $|M| = m$ implies that $M$ contains a circuit, so we can apply Lemma~\ref{lemStructuralCircuits}.
Suppose $m=3$. 
The following claim shows that \eqref{eqn4ptlineplus3} holds.
\begin{claim}\label{claim_m=3}
$\bA\sim \left[\BE^1|\BE^1+2\BE^2|\BE^2|\BE^1+\BE^2|\BE^3|\BE^1+\BE^3|\BE^1+2\BE^2+\BE^3|2\BE^1+2\BE^2+\BE^3|\BE^1+\BE^2+\BE^3\right]$.
\end{claim}
\proof{Proof of Claim.}
By Lemma~\ref{lemStructuralCircuits}~\emph{\ref{struc1}} and~\eqref{eqn:C0}, the circuit has the form $ [\BE^1_{2}|\bb|\bb+2\BE^1_{2}]$ for some $\bb \in \bA/ \BE^1$.
After elementary operations, we have
\[
\bA \supseteq \left[
\begin{array}{c@{\hskip .1 cm}|@{\hskip .1 cm}c@{\hskip .1 cm}|@{\hskip .1 cm}c@{\hskip .1 cm}|@{\hskip .1 cm}c@{\hskip .1 cm}|@{\hskip .1 cm}c@{\hskip .1 cm}|@{\hskip .1 cm}c@{\hskip .1 cm}|@{\hskip .1 cm}c}
 \multirow{2}{*}{$\bB^*$}& \beta_0 & \beta_0+1 & \beta_1 & \beta_1+1 & \beta_2 & \beta_2+1  \\
 & \BE^1_2 & \BE^1_2 & \bb & \bb & \bb+2\BE^1_{2} & \bb+2\BE^1_{2}
\end{array}
\right] \sim
\left[
\begin{array}{c@{\hskip .1 cm}|@{\hskip .1 cm}c@{\hskip .1 cm}|@{\hskip .1 cm}c@{\hskip .1 cm}|@{\hskip .1 cm}c@{\hskip .1 cm}|@{\hskip .1 cm}c@{\hskip .1 cm}|@{\hskip .1 cm}c@{\hskip .1 cm}|@{\hskip .1 cm}c@{\hskip .1 cm}|@{\hskip .1 cm}c}
1 & 1 & \beta_0 & \beta_0+1 &\beta_1&\beta_1+1&\beta_2&\beta_2+1\\
0 & 2 & 1 & 1 &0&0 & 2 & 2\\
0 & 0 & 0 & 0 &1&1 & 1 & 1
\end{array}
\right].
\]
We assume, without loss of generality, that the equivalence in the latter displayed equation is an equation.
We must have $\beta_0 =0 $ because it is the midpoint of the columns of $\bB^*$; see Lemma~\ref{lem:BasicProperties}~\emph{\ref{propBasicProperties3}}.
By subtracting the third row from the first row $\beta_1$ many times, we assume $\beta_1 = 0$; we then see that $2|\beta_2-2| \le 2$ using bimodularity with columns 2, 6, and 7 and that $2|\beta_2|\le 2$ using bimodularity with columns 2, 5, and 8.
Hence, $\beta_2 = 1$. 
By Lemma~\ref{lem:BasicProperties}~\emph{\ref{propBasicProperties3}} and the maximality of $\bA$, it follows that $\sfrac{1}{2}\cdot\left( (\beta_2+1,2,1)+(\beta_1,0,1)\right) = (1,1,1) \in \bA$.
Thus, $\bA$ has at least $\cc(2,3)= 9$ columns. 

The following equation follows from the definition of $\bA/ \BE^1$:
\[
|\bA| = |\bO(\bzero)| + |\bA/ \BE^1| + \sum_{\bb \in \bA/ \BE^1}\big ( | \bO(\bb) | - 1 \big).
\]
By using~\eqref{M=2} and Lemma~\ref{lem:BasicProperties}~\emph{\ref{propBasicProperties1}}, we see that
\[
|\bA|=|\bO(\bzero)| + |\bA/ \BE^1| + \sum_{\bb \in \bA/ \BE^1}\big ( | \bO(\bb) | - 1 \big) = 1 + |\bA/ \BE^1| + |\bM|.
\]
The matrix $\bA/ \BE^1 \subseteq \Z^2$ is bimodular and contains $2\BE^1_2$ and $\BE^1_2$. 
It is quickly verified that (after multiplying columns by $-1$), we have $|\bA/ \BE^1|\le 5$;
%
%
see also the proof of Proposition~\ref{prop:m2} for an argument of this. 
Hence, 
\begin{equation}\label{eqm3}
|\bA| = 1+|\bA/ \BE^1|+|\bM| \le 1+5+3 = 9.
\end{equation}
Thus, $|\bA| = 9$ and $\bA$ can be transformed via elementary operations to the form in Claim~\ref{claim_m=3}.
\hfill$\diamond$
\endproof

\smallskip
Suppose $m \ge 4$.
%
%
%
%
It follows from Lemma~\ref{lemStructuralCircuits}~\emph{\ref{struc1}} that $\bM = [\BE^1_{m-1}|\bb^1|\cdots|\bb^{m-1}]$ and $\bM$ contains exactly one circuit $[\BE^1_{m-1}|\bb^1|\bb^2] =[\BE^1_{m-1}|\bb^1|\bb^1+2\BE^1_{m-1}]$.
Hence,
\[
\bE := 
\left[
\begin{array}{c@{\hskip.1cm}|@{\hskip.1cm}c@{\hskip.1cm}|@{\hskip.1cm}c@{\hskip.1cm}|@{\hskip.1cm}c}
\bb^1&\bb^2&\cdots&\bb^{m-1}
\end{array}
\right]
\]
is a basis and $|\det \bE|= |\det[\bb^1|\bb^2|\cdots|\bb^{m-1}]| = |\det[\bb^1|2\BE^1_{m-1}|\cdots|\bb^{m-1}]| =  2$.
Observe that 
\[
\bE\sim 
 \left[
 \begin{array}{c@{\hskip.1cm}|@{\hskip.1cm}c@{\hskip.1cm}|@{\hskip.1cm}c@{\hskip.1cm}|@{\hskip.1cm}c@{\hskip.1cm}|@{\hskip.1cm}c}
 \BE^2_{m-1}&\BE^2_{m-1}+2\BE^1_{m-1}&\BE^3_{m-1}&\cdots&\BE^{m-1}_{m-1}
 \end{array}\right],
\]
so if $\bE \bw \in \bA/ \BE^1$ for some $\bw  = (w_1, \dotsc, w_{m-1})\in \R^{m-1}$, then $w_3, \dotsc, w_{m-1} \in \{-1,0,1\}$.
 
For $i = 1, \dotsc, m-1$, let $\beta_i\in \Z$ be such that $(\beta_i, \bb^i), (\beta_i+1, \bb^i) \in \bA$.
Define
\[
\Gamma := \{\bgamma = (\gamma_{1}, \dotsc, \gamma_{m-1}) \colon \gamma_i \in \{\beta_i, \beta_i+1\} ~\forall ~ i = 1, \dotsc, m-1\}.
\]
For each $(\beta, \bb) \in \bA$ and $\bgamma = (\gamma_1,\dotsc, \gamma_{m-1}) \in \Gamma$, Cramer's rule and bimodularity of $\bA/\BE^1$ imply $\bE^{-1}\bb =: \bv = (v_1, \dotsc, v_{m-1})\in\{ \pm 1, \pm \sfrac{1}{2} , 0\}^{m-1}$.
Furthermore,
\[
\left| \det 
\left[ 
\begin{array}{c@{\hskip .1 cm}|@{\hskip .1 cm}c} 
 \bgamma^\intercal &\beta \\
 \bE & \bb
\end{array}\right]
\right| 
=
\left| \det 
\left[ 
\begin{array}{c@{\hskip .1 cm}|@{\hskip .1 cm}c} 
 \bgamma^\intercal &\beta - \bgamma^\intercal \bv \\
 \bE & \bzero
\end{array}\right]
\right| 
=
2  |\beta - \bgamma^\intercal \bv|
\le 2,
\]
so $|\beta - \bgamma^\intercal \bv| \le 1$. 
As $v_3, \dotsc, v_{m-1}\in \{-1,0,1\}$, the maximum $\sigma_{\max}$ and minimum $\sigma_{\min}$ of $\{\beta - \bgamma^\intercal \bv \colon \bgamma \in \Gamma\}$ satisfy $2\ge\sigma_{\max} - \sigma_{\min} \ge \sfrac12\cdot |\supp(v_1,v_2)| + |\supp (v_3, \dotsc, v_{m-1})|$. 
Thus, we have $|\supp (v_3, \dotsc, v_{m-1})| \le 2$ because $|\beta - \bgamma^\intercal \bv| \le 1$ holds for all $\bgamma \in \Gamma$.
This leads to a natural partition of $\bA$.
For $j = 0,1,2$, define
\[
\bA^j  :=  \left\{ (\beta, \bb) \in \bA \colon  
 \bE^{-1}\bb = (v_1, \dotsc, v_{m-1})~\text{with}~|\supp (v_3, \dotsc, v_{m-1})| = j  \right\}.
\]
We have $\bM \subseteq [\bA^0|\bA^1]$ because $\bM = [\BE^1_{m-1}|\bb^1|\cdots|\bb^{m-1}]$ and $[\BE^1_{m-1}|\bb^1|\bb^2]$ is a circuit.

\begin{claim}\label{claimA0}
$\bA^0\sim [\BE^1|\BE^1+2\BE^2|\BE^2|\BE^1+\BE^2|\BE^3|\BE^1+\BE^3|\BE^1+2\BE^2+\BE^3|2\BE^1+2\BE^2+\BE^3|\BE^1+\BE^2+\BE^3]$.
\end{claim}
\proof{Proof of Claim.}
By definition, $\bA^0 = \bA \cap \Span\{\BE^1, (\beta_1, \bb^1), (\beta_2, \bb^1+2\BE^1_{m-1})\}$ and $\rank \bA^0 = 3$.
Recall~\eqref{eqProject}: $\bA^0 \sim (\overline{\bA^0}, \bzero)$, where $\overline{\bA^0}$ is a full row rank bimodular matrix. 
One such sequence of elementary operations maps $\BE^1$ to $\BE^1_{3}$, $\BE^2$ to $\BE^2_3$, and $\bb^1$ to $\BE^3_{3}$.
By Claim \ref{claim_m=3}, we know that $ |\bA^0| = 9$ and $\bA^0$ can be transformed via elementary operations to the form described in Claim~\ref{claimA0}.
%
\hfill$\diamond$
\endproof

\begin{claim}\label{claimA1}
For each $i = 3, \dotsc, m-1$, there are at most four columns $(\beta, \bb) \in \bA$ such that $\bE^{-1}\bb =  (v_1, \dotsc, v_{m-1})$ satisfies $\supp (v_3, \dotsc, v_{m-1}) = \{i\}$.
Consequently, $|\bA^1| \le 4(m-3)$. 
\end{claim}
\proof{Proof of Claim.}
Following Claim~\ref{claimA0}, we assume
\begin{equation}\label{eqA0transfrom}
\bA^0= [\BE^1|\BE^1+2\BE^2|\BE^2|\BE^1+\BE^2|\BE^3|\BE^1+\BE^3|\BE^1+2\BE^2+\BE^3|2\BE^1+2\BE^2+\BE^3|2\BE^1+\BE^2+\BE^3].
\end{equation}
Set 
\[
\bF^i := \left\{(\beta, \bb) \in \bA \colon  \bE^{-1}\bb =  (v_1, \dotsc, v_{m-1})~\text{satisfies}~\supp (v_3, \dotsc, v_{m-1}) = \{i\}\right\}.
\]
We claim that
\begin{equation}\label{eqDifferenceA1}
\text{$\ba - \ba'\in [\bA^0|-\bA^0]$ for every pair of distinct columns $\ba, \ba' \in \bF^i$.}
\end{equation}
Assume to the contrary that~\eqref{eqDifferenceA1} is violated by some $\ba, \ba' \in F^i$.
Recall~\eqref{eqProject}: $[A^0|\ba-\ba'] \sim (\overline{\bF},\bzero)$, where $\overline{\bF}$ has full row rank and differing columns.
We have $\rank \overline{\bF} =3$ because $\ba, \ba' \in F^i$.
Claim~\ref{claim_m=3} established that a rank-3 bimodular matrix $\bA$ containing $\bA^0$ has at most nine differing columns.
It follows that $\overline{\bF}$ is not bimodular.
%
%
In particular, there exists a basis $\left[\overline{\ba\vphantom{\bd}}-\overline{\ba'\vphantom{\bd}}|\overline{\bc\vphantom{\bd}}|\overline{\bd}\right]\subseteq \overline{\bF}$ such that $|\det\left[\overline{\ba\vphantom{\bd}}-\overline{\ba'\vphantom{\bd}}|\overline{\bc\vphantom{\bd}}|\overline{\bd}\right]| \ge 3 $.
By~\eqref{eqPreserve}, any basis in $\bA$ containing $[\ba|\ba'|\bc|\bd]$ has an absolute determinant of at least three, which contradicts that $\bA$ is bimodular.
This shows that~\eqref{eqDifferenceA1} is true.

Note that $\bF^i$ contains $(\beta_i, \bb^i)$ and $(\beta_i+1, \bb^i)$.
Let $\ba, \ba' \in \bF^i \setminus \{(\beta_i, \bb^i),(\beta_i+1, \bb^i)\}$.
The following are columns of $[\bA^0|-\bA^0]$ according to~\eqref{eqDifferenceA1}: $\ba - (\beta_i, \bb^i)$, $\ba - (\beta_i+1,\bb^i) = \ba - (\beta_i,\bb^i) - \BE^1$, $\ba'- (\beta_i, \bb^i)$, and $\ba' - (\beta_i+1,\bb^i) = \ba' - (\beta_i,\bb^i) - \BE^1$.
From~\eqref{eqA0transfrom} and the assumption that $A$ only contains primitive columns, it follows that
\[
    \ba-\left[\begin{array}{c}\beta_i\\\bb^i\end{array}\right], ~\ba'-\left[\begin{array}{c}\beta_i\\\bb^i\end{array}\right] \in\left\{\BE^1+\BE^2,\BE^1+\BE^3,-(\BE^1+2\BE^2+\BE^3),2\BE^1+2\BE^2+\BE^3,-\BE^2,-\BE^3\right\}.
\]
Furthermore, $\ba - \ba' = \left(\ba - (\beta_i, \bb^i)\right) - \left(\ba' - (\beta_i, \bb^i)\right)$ is a column of $[\bA^0|-\bA^0]$.
Define $\bS^1 := \left\{\BE^1+\BE^2,\BE^1+\BE^3,-(\BE^1+2\BE^2+\BE^3)\right\}$ and $\bS^2 := \left\{2\BE^1+2\BE^2+\BE^3,-\BE^2,-\BE^3\right\}$.
It is quickly checked that if both $\left(\ba - (\beta_i, \bb^i)\right) $ and $\left(\ba' - (\beta_i, \bb^i)\right)$ are in $\bS^1$ or both are in $\bS^2$, then $\ba - \ba' \not \in [\bA^0|-\bA^0]$.
Hence, there are at most two columns $\bF^i \setminus \{(\beta_i, \bb^i),(\beta_i+1, \bb^i)\}$.
Equivalently, $|\bF^i| \le 4$.
\hfill$\diamond$
\endproof

\begin{claim}\label{claimA2}
$\left|\bA^2\right| \le \binom{m-3}{2}$.
\end{claim}
\proof{Proof of Claim.}
Let $(\beta, \bb) \in \bA^2$.
Set $\bE^{-1}\bb =: \bv = (v_1, \dotsc, v_{m-1})$, where $\supp(v_3, \dotsc, v_{m-1}) = \{i, j\}$ for some $i,j \in \{3, \dotsc, m-1\}$.
To prove the claim, it suffices to show that there is no other column $(\beta', \bb') \in \bA^2$ such that $\bE^{-1}\bb' =: \bv ' =  (v'_1, \dotsc, v'_{m-1})$ satisfies $\supp \bv' = \{i, j\}$. 
Indeed, this will show that a column of $\bA^2$ is uniquely determined by two indices in $\{3, \dotsc, m-1\}$.

For simplicity, assume $i=3$ and $j=4$.
 Recall $2\ge \sigma_{\max} -\sigma_{\min} \ge \sfrac12\cdot |\supp(v_1,v_2)| + |\supp(v_3,\dotsc,v_{m-1})|\ge|\supp(v_3,v_4)|= 2$.
Therefore, it must hold that $\supp \bv = \{3, 4\}$; in particular, $v_1 = v_2 = 0$. 
Assume to the contrary that $\bA$ contains another column $(\beta', \bb')$ such that $\bE^{-1}\bb' =: \bv ' =  (v'_1, \dotsc, v'_{m-1})$ satisfies $\supp \bv' = \{3,4\}$.
As was the case with $\bv$, we have $\supp \bv' = \{3,4\}$.
Recall that $\bM \subseteq [\bA^0|\bA^1]$. 
This implies $\bb \not \in \bM$ and $|\bO(\bb)| = 1$.
Therefore, $\bb$ and $\bb'$ are distinct because $|\bO(\bb')| = 1$ and $(\beta, \bb)\neq (\beta', \bb')$.
In fact, the two columns differ according to assumption~\eqref{eqAssumeNoNegatives}.
Given that $v_3, v_4, v'_3, v'_4 \in \{-1,1\}$, the columns $(v_3, v_4)$ and $(v'_3, v'_4)$ must have different sign patterns, say $v_3 = v_3' = 1$ and $v_4 = -v'_4 = 1$. 
By multi-linearity of the determinant,
\[
\left|\det\left[
\begin{array}{c@{\hskip.1cm}|@{\hskip.1cm}c@{\hskip.1cm}|@{\hskip.1cm}c@{\hskip.1cm}|@{\hskip.1cm}c@{\hskip.1cm}|@{\hskip.1cm}c@{\hskip.1cm}|@{\hskip.1cm}c@{\hskip.1cm}|@{\hskip.1cm}c}
\bb^1&\bb^2&\bb&\bb'&\bb^5&\cdots&\bb^{m-1}
\end{array}
\right]\right| 
= \left|\det\left[
\begin{array}{c@{\hskip .1 cm}|@{\hskip .1cm}c}
v_3 & v'_3\\
v_4 & v'_4
\end{array}
\right]
\right|
\cdot \left|\det \left[
\begin{array}{c@{\hskip.1cm}|@{\hskip.1cm}c@{\hskip.1cm}|@{\hskip.1cm}c}
\bb^1&\cdots&\bb^{m-1}
\end{array}
\right]\right |
= 2  |\det E| = 4.
\]
This contradicts that $\bA/ \BE^1$ is bimodular. 
\hfill$\diamond$
\endproof

\smallskip
Finally, we prove~\eqref{eqn4ptlineplus3} by combining Claims~\ref{claimA0},~\ref{claimA1}, and~\ref{claimA2}:
    \[
    |\bA| = |\bA^0| + |\bA^1| + |\bA^2| 
    \le 9 + 4(m-3) + \binom{m-3}{2}
    =  \frac{1}{2} (m^2+m) + 3.
    \]
\hfill\halmos
\endproof

\subsection{Additional structural properties when \texorpdfstring{$|\bB^*| = 3$}~.}\label{subsecExtraProperties3}

Throughout Section~\ref{secCircuitStructure}, we have assumed that $A$ has the form~\eqref{eqDefnA}. 
In order to transform $A$ into this form, we identify a primitive column to transform into $\BE^1$.
Our choice of a primitive column thus far has been somewhat arbitrary when in reality there are multiple choices.
For example, if $[\BE^1|\BE^2|\BE^1+\BE^2+2\BE^3] \subseteq \bA$, then any of these columns can be chosen.
Moreover, these columns are interchangeable: if we label the columns as $\ba^1, \ba^2, $ and $\ba^3$, then for any permutation $\sigma \in \mathcal{S}^3$ there are elementary operations such that $[\ba^{\sigma(1)}|\ba^{\sigma(2)}|\ba^{\sigma(3)}] = [\BE^1|\BE^2|\BE^1+\BE^2+2\BE^3]$.
The discussion surrounding Equation~\eqref{eqA2circuits} in the proof of Lemma~\ref{leme1e2e3Gen} illustrates this symmetry in more detail.
In Lemma~\ref{leme1e2e3Gen}, we consider swapping the roles of $\BE^1$ and another primitive column in $[\BE^1|\BE^2|\BE^1+\BE^2+2\BE^3]$.
In order to formalize our argument, we define $\bA/\ba$ for a primitive column $\ba \in \bA$ to be the matrix $\bA/\BE^1$ after identifying $\ba$ with $\BE^1$.
We use $\bM_{\ba}$ to denote the set of columns of $\bA/ \ba$ with at least two originals in $\bA$, and we use $\bO(\bC)$ to denote the set of original columns in $\bA$ corresponding to a subset $\bC \subseteq \bA/ \ba$.

\begin{lemma}\label{leme1e2e3Gen}
Assume that $\bA$ only contains primitive columns and $\bM$ contains a circuit.
Choose $\bB^*\subseteq \bA$ that satisfies~\eqref{eqMinimalAssumption} and minimizes $|\bB^*|$, and assume that $\bB^*$ has the form~\eqref{eqBForm}.
If $|\bB^*| = 3$, then there exists at least one column $\ba \in [\BE^1|\BE^2|\BE^1+\BE^2+2\BE^3]$ such that $\bM_{\ba}$ does not contain a circuit. 
It follows that $|M_\ba| \le m-1$.
\end{lemma}

\proof{Proof.}
Let $\ba \in [\BE^1|\BE^2|\BE^1+\BE^2+2\BE^3]$.
As stated in the previous paragraph, we can assume $\ba = \BE^1$ by applying elementary operations to $\bA$.
Lemma~\ref{lemStructuralCircuits}~\emph{\ref{struc2}} implies $|\bM_{\ba}| \le m$ for each $\ba \in [\BE^1|\BE^2|\BE^1+\BE^2+2\BE^3]$.

\begin{claim}\label{prop2:leme1e2e3Gen}
Suppose there exists a column $\ba \in [\BE^1|\BE^2|\BE^1+\BE^2+2\BE^3]$ such that $\bM_{\ba}$ contains a circuit $\bC$; then there exists a column $\ba'\in [\BE^1|\BE^2|\BE^1+\BE^2+2\BE^3] \setminus \{ \ba\}$ such that $\bM_{\ba'}$ contains a circuit $\bC'$.
Furthermore, $\bO(\bC) = \bO(\bC')$.
\end{claim}
\proof{Proof of Claim.}
Without loss of generality, $\ba = \BE^1$.
It follows from Lemma~\ref{lemStructuralCircuits}~\emph{\ref{struc2}} that $\bC = [ \BE^1_{m-1}+\BE^2_{m-1}|\BE^2_{m-1}|\BE^3_{m-1}|\bd+\BE^3_{m-1}]$ for some $\bd = (d_1, d_2, \bzero) \in[\BE^1_{m-1}|\BE^1_{m-1}+2\BE^2_{m-1}]$.
Given that $\bC \subseteq \bM$, we have $(\alpha, \BE^2_{m-1}+\BE^3_{m-1})$ and $(\alpha+1, \BE^2_{m-1}+\BE^3_{m-1})$ are in $\bA$ for some $\alpha \in \Z$.
By~\eqref{BSpan}, we can assume $\alpha = 0$, so $\BE^2+\BE^3, \BE^1+\BE^2+\BE^3 \in \bA$. 
Similarly, we can assume $\BE^3,\BE^1+\BE^3\in \bA$.
Suppose $(\beta, \bd+\BE^3_{m-1}), (\beta+1, \bd+\BE^3_{m-1}) \in \bA$.
Hence, 
\begin{equation}\label{eqA2circuits}
A 
= [\bB^*|~\bO(\bC)~|~A']=
\left[
\begin{array}{c@{\hskip .1 cm}|@{\hskip .1 cm}c}
\begin{array}{c@{\hskip .1 cm}|@{\hskip .1 cm}c@{\hskip .1 cm}|@{\hskip .1 cm}c@{\hskip .1 cm}|@{\hskip .1 cm}c@{\hskip .1 cm}|@{\hskip .1 cm}c@{\hskip .1 cm}|@{\hskip .1 cm}c@{\hskip .1 cm}|@{\hskip .1 cm}c@{\hskip .1 cm}|@{\hskip .1 cm}c@{\hskip .1 cm}|@{\hskip .1 cm}c@{\hskip .1 cm}|@{\hskip .1 cm}c@{\hskip .1 cm}|@{\hskip .1 cm}c}
1 & 0 & 1 & 0 & 1  & 0 & 1 & 0 & 1 & \beta & \beta+1  \\
0 & 1 & 1 & 1 & 1  & 0 & 0 & 0 & 0 & d_1 & d_1 \\
0 & 0 & 2 & 1 & 1  & 1 & 1 & 0 & 0 & d_2 & d_2 \\
0 & 0 & 0 & 0 & 0  & 0 & 0 & 1 & 1 & 1 & 1 \\
\bzero &\bzero &\bzero & \bzero &  \bzero  & \bzero  & \bzero & \bzero   & \bzero & \bzero & \bzero 
\end{array}
&
A'
\end{array}
\right].
\end{equation}
Columns 4-11 in~\eqref{eqA2circuits} correspond to $\bO(\bC)$.

Suppose $\bd = \BE^1_{m-1}$, or equivalently suppose $(d_1, d_2) = (1,0)$.
It follows that $\beta=0$, and Columns 4-11 correspond to $\bO(\bC')$, where $\bC' \subseteq \bM_{\BE^2}$ is a circuit. 
This proves the claim.

Suppose $\bd=\BE^1_{m-1}+2\BE^2_{m-1}$ or equivalently suppose $(d_1, d_2) = (1,2)$.
From~\eqref{eqPreserve}, the top left $4\times 11$ submatrix of~\eqref{eqA2circuits} is bimodular.
From this, we see that $\beta  = 1$: if $\beta \le 0$, then the first four rows of columns 2,7,9, and 10 form a basis with absolute determinant greater than two, and if $\beta \ge 2$, then the first four rows of columns 2,6,8, and 11 form a basis with absolute determinant greater than two.
Conduct the following three elementary row operations to $\bA$ followed by multiplying columns by $-1$: (1) subtract the second row from the first row and subtract twice the second row from the third row; (2) multiply the third row by $ -1$; (3) add the third row to the first row; (4) negate columns 6 and 7:
\begin{equation}\label{eq2ndapex}
A 
\sim
\left[
\begin{array}{c@{\hskip .1 cm}|@{\hskip .1 cm}c}
\begin{array}{c@{\hskip .1 cm}|@{\hskip .1 cm}c@{\hskip .1 cm}|@{\hskip .1 cm}c@{\hskip .1 cm}|@{\hskip .1 cm}c@{\hskip .1 cm}|@{\hskip .1 cm}c@{\hskip .1 cm}|@{\hskip .1 cm}c@{\hskip .1 cm}|@{\hskip .1 cm}c@{\hskip .1 cm}|@{\hskip .1 cm}c@{\hskip .1 cm}|@{\hskip .1 cm}c@{\hskip .1 cm}|@{\hskip .1 cm}c@{\hskip .1 cm}|@{\hskip .1 cm}c}
1 & 1 & 0 & 0 & 1  & 1 & 0 & 0 & 1 & 0 & 1  \\
0 & 1 & 1 & 1 & 1  & 0 & 0 & 0 & 0 & 1 & 1 \\
0 & 2 & 0 & 1 & 1  & 1 & 1 & 0 & 0 & 0 & 0 \\
0 & 0 & 0 & 0 & 0  & 0 & 0 & 1 & 1 & 1 & 1 \\
\bzero &\bzero &\bzero & \bzero &  \bzero  & \bzero  & \bzero & \bzero   & \bzero & \bzero & \bzero 
\end{array}
&
A''
\end{array}
\right].
\end{equation}
Columns 4-11 in~\eqref{eqA2circuits} correspond to $\bO(\bC)$, where $\bC \subseteq \bA / \BE^1$ is a circuit.
Similarly, Columns 4-11 in~\eqref{eq2ndapex} correspond to $\bO(\bC')$, where $\bC' \subseteq \bA / \BE^2$ is a circuit.
Furthermore, Columns 4-11 in~\eqref{eqA2circuits} are equivalent up to row and column operations to Columns 4-11 in~\eqref{eq2ndapex}.
Now, $\BE^2$ on the right hand side of~\eqref{eq2ndapex} is equivalent to $\BE^1+\BE^2+2\BE^3$ in~\eqref{eqA2circuits}. 
Therefore, in the original representation of $\bA$ in~\eqref{eqA2circuits}, we conclude that there is a circuit $\bC' \subseteq \bM_{ \BE^1+\BE^2+2\BE^3}$ such that $\bO(\bC') = \bO(\bC)$.
\hfill $\diamond$
\endproof

Assume to the contrary that $\bM_{\BE^1}$, $\bM_{\BE^2}$, and $\bM_{\BE^1+\BE^2+2\BE^3}$ each contain a circuit.
Call these circuits $\bC^1, \bC^2$, and $\bC^3$, respectively. 
By Claim~\ref{prop2:leme1e2e3Gen}, for each $i \in \{ 1, 2,3\}$, there exists some $j_i \in\{1,2,3\}\setminus \{i\}$ such that $\bO(\bC^i) = \bO(\bC^{j_i})$.
Since there is an odd number of circuits here, we conclude that $\bO(\bC^1)=\bO(\bC^2) =\bO(\bC^3)$.
This means that columns 4-11 in~\eqref{eqA2circuits} equal $\bO(\bC^1),\bO(\bC^2),$ and $\bO(\bC^3)$.
Suppose $(d_1, d_2) = (1,0)$.
Denote column 4 by $\bc$. 
There is no column $\bc'$ among 5-11 such that $\bc - \bc' \in \{ \pm (\BE^1+\BE^2+2\BE^3)\}$.
Therefore, $\bc$ is not an original column of $\bC^3$, which is a contradiction. 
Similarly, if $(d_1, d_2) = (1,2)$, then we see that column $10$ is not an original column of $\bO(\bC^2)$.
This proves Lemma~\ref{leme1e2e3Gen}.
\hfill\halmos
\endproof

\medskip

This ends our discussion on new combinatorial properties of bimodular constraint matrices. 
We reiterate to the reader that we believe Lemmas~\ref{lemStructuralCircuits} and~\ref{leme1e2e3Gen} may be of independent interest in future research.
Next, we apply these properties to prove Theorem~\ref{thm:bimodular}.

\section{A proof of Theorem~\ref{thm:bimodular}.}\label{secDelta=2}

Proposition~\ref{prop:lowerbounds} proves $\cc(2,m) \ge \sfrac{1}{2} \cdot (m^2+m)+m$.
We prove $\cc(2, m) = \sfrac{1}{2} \cdot \left(m^2+m\right)+m$ using induction on $m$.
It is quickly verified that when $m=1$ the unique maximal bimodular matrix is $[1|2]$ (up to multiplying columns by $-1$).
This proves $\cc(2,1) = 2$.
Assume that $m \ge 2$ and
\begin{equation}\label{eqInductiveStep}
\cc(2, k) = \frac{1}{2}  \left(k^2+k\right)+k \qquad \forall  ~k = 1, \dotsc, m-1.
\end{equation}
Let $A \in \Z^{m\times n}$ be a maximal bimodular matrix with $\rank \bA = m$ and differing columns.
After elementary operations, we can assume that $\BE^1 \in \bA$.
For the inductive step, we use the following relationship between $|\bA|$ and $|\bA/ \BE^1|$:
\[
|\bA| = |\bO(\bzero)| + |\bA/ \BE^1| + \sum_{\bb \in \bA/ \BE^1}\big ( | \bO(\bb) | - 1 \big).
\]
Recall $\bM = \{\bb \in \bA/ \BE^1 \colon |\bO(\bb)| \ge 2\}$. 
We consider two cases: when $\bA$ only contains primitive columns and when $\bA$ contains a non-primitive column.

\smallskip
\noindent{\bf Case 1.}
Assume that $\bA$ only contains primitive columns. 
We have $|\bO(\bzero)|=1$ by Lemma~\ref{lem:BasicProperties}~\emph{\ref{propBasicProperties1}}.

\smallskip
\noindent{\bf Subcase 1.1.}
Assume that $\bM$ contains a circuit or $\bM$ contains a column $\bb'$ with $|\bO(\bb')|=3$. 
If $|\bO(\bb')|=3$, then by Lemma~\ref{lem:BasicProperties}~\emph{\ref{propBasicProperties3}} we have $B^*\subseteq A$ that satisfies~\eqref{eqMinimalAssumption} and $|B^*|=2$. If $M$ contains a circuit, then we
choose $\bB^*\subseteq \bA$ that satisfies~\eqref{eqMinimalAssumption} and minimizes $|\bB^*|$. 
We have $2 \le |\bB^*|\le 4$ according to Lemma~\ref{lemSmallCircuits}.
In both cases, we can assume that $\bB^*$ has the form~\eqref{eqBForm}, and $|\bM| = \sum_{\bb\in \bA/ \BE^1} (|\bO(\bb)|-1)$ by \eqref{M=2}.

Suppose $|\bB^*| = 2$.
If $m = 2$, then $A \subseteq \Span \bB^*$ by Lemma~\ref{lemSpan}.
It follows from~\eqref{BSpan} that $A \subseteq \bA \cap \Span \bB^* \subseteq [\bB^* | \BE^2| \BE^1+\BE^2]$.
Hence, $|\bA| \le 4 = \sfrac{1}{2} \cdot (m^2+m)+m-1$.
Suppose $m \ge 3$.
Lemma~\ref{lemStructuralCircuits}~\emph{\ref{struc1}} implies $|\bM| \le m$.
If $|\bM|=m$, then $|\bA|\le \sfrac{1}{2}\cdot  (m^2+m) +3$ by Lemma~\ref{lem4ptlineplus3}.
If $|\bM|\le m-1$, then 
\[
|\bA|=1+|\bA/\BE^1|+|\bM|\le 1+\cc(2,m-1)+(m-1) \le \sfrac{1}{2}\cdot \left(m^2+m\right)+m-1.
\]
Therefore,
\begin{equation}\label{eqe1e2}
|\bA|\le
\begin{cases}
        \sfrac{1}{2}\cdot (m^2+m)+m-1,& \text{if}~m\ge 4~\text{or}~ m=2;\\[.1 cm]
        \sfrac{1}{2}\cdot(m^2+m)+m, & \text{if}~m=3.
    \end{cases}
\end{equation}

Suppose $|\bB^*| =3$.
We have $[\BE^1_{m-1}|\BE^1_{m-1}+2\BE^2_{m-1}]\subseteq \bA/\BE^1$ and every column of $\bA/\BE^1$ is primitive because $2 < |\bB^*| = 3$. 
Therefore, $\bA/\BE^1$ satisfies Inequality~\eqref{eqe1e2}.
By Lemma~\ref{leme1e2e3Gen}, there exists a column $\ba \in[\BE^1| \BE^2| \BE^1+\BE^2+2\BE^3]$ such that $|\bM_{\ba}| = \sum_{\bb\in \bA/\ba}(|\bO(\bb)|-1) \le m-1$.
Therefore,
\begin{equation}\label{eqe1e2e3}
%
|\bA| \le
 \begin{cases}
        \sfrac{1}{2}\cdot (m^2+m) +m-2, & \text{if}~m\ge 5~\text{or}~ m=3;\\
        \sfrac{1}{2}\cdot (m^2+m)+m-1, &\text{if}~ m=4.
    \end{cases}
\end{equation}

Suppose $|\bB^*| = 4$.
If $m=4$, then $|\bA|\le 12 = \sfrac{1}{2}\cdot(m^2+m)+m-2$ by $\eqref{BSpan}$.
Suppose $m\ge 5$; then $|\bM|\le m+1$ by Lemma~\ref{lemStructuralCircuits}~\emph{\ref{struc3}}.
 The matrix $\bB^*/\BE^1$ satisfies $|\bB^* / \BE^1| = 3$. 
Furthermore, $\bB^* / \BE^1$ has minimal cardinality among all subsets of $\bA / \BE^1$ satisfying~\eqref{eqMinimalAssumption} otherwise we would contradict the minimality of $B^*$.
Hence, $\bA/ \BE^1$ satisfies Inequality~\eqref{eqe1e2e3} for $m-1$.
Therefore,
\begin{equation}\label{eqe1e2e3e4}
|\bA| \le \begin{cases}
        \sfrac{1}{2}\cdot(m^2+m)+m-1, & \text{if}~m\ge 6;\\
        \sfrac{1}{2}\cdot(m^2+m)+m, & \text{if}~m=5;\\
        \sfrac{1}{2}\cdot(m^2+m)+m-2, & \text{if}~m=4.
    \end{cases}
\end{equation}

\noindent{\bf Subcase 1.2.} 
Assume that $\bM$ does not contain a circuit and $\bM = \{\bb \in \bA/ \BE^1 \colon |\bO(\bb)|=2\}$. 
%
This implies $|\bM| \le m-1$ and $|\bA|  \le 1+\cc(2,m-1)+(m-1) \le \sfrac{1}{2}\cdot \left(m^2+m\right)+m-1$.
\smallskip
\noindent{\bf Case 2.} Assume that $\bA$ contains a non-primitive column $\ba$.
The column $\sfrac{1}{2}\cdot \ba$ is contained in $A$ because $\bA$ is maximal, and the column is primitive by Lemma~\ref{lem:BasicProperties}.
By transforming $\sfrac{1}{2} \cdot \ba$ to $\BE^1$ using elementary operations, $\ba$ transforms to $2\BE^1$, and we can write
\[
\bA =
\left[\begin{array}{c@{\hskip .15cm}|@{\hskip .15cm}c}
    2\BE^1 & \bA'
\end{array}
\right]=
\left[
\begin{array}{c@{\hskip .1 cm}|@{\hskip .1 cm}c}
2 & \bbeta^\intercal \\
\bzero & \hat{\bA}
\end{array}
\right],
\]
where $\bA' \in \Z^{m \times (n-1)}$ and $\BE^1 \in \bA'$.
From this identity we see that $\hat{\bA} \supseteq \bA/ \BE^1$ is unimodular, so $|\bA/ \BE^1| \le \cc(1,m-1) = \sfrac{1}{2}\cdot  (m^2-m)$.
We refer to known results in matroid theory to complete the case; see~\cite{Oxley1992} for a thorough introduction on matroids.
From Lemma~\ref{lem:BasicProperties}~\emph{\ref{propBasicProperties1}} and~\emph{\ref{propBasicProperties3}}, it follows that $\ba - \ba' \not\in p \cdot \Z^m$ for any distinct columns $\ba, \ba' \in \bA'$ and any prime number $p \ge 3$.
This, along with the assumption that $\bA'$ is bimodular, demonstrates that the matrix $\bA'$ is a representation of a matroid $\cM$ over the field $\text{GF}(p)$ for any prime number $p \ge 3$.
Similarly, the matrix $\bA/ \BE^1$ is a representation of the simplification $\cM / \BE^1$ of the minor of $\cM$ obtained by contracting $\BE^1$; here we use the fact that $\bA/ \BE^1$ is defined to have differing columns. 
According to Kung~\cite[Lemma 2.2.1]{KunMat1990}, $\cM$ does not contain the Reid geometry.
By~\cite[Theorem 3.1]{KunMat1990}, it follows that $|\bA'|  -|\bA/ \BE^1| =  |\cM|-|\cM/\BE^1| \le 2m-1$ .
Therefore, 
\[
|\bA| = 1+\left|\bA'\right| \le 1 + (2m-1)+|\bA/\BE^1| \le 2m + \frac{1}{2}(m^2-m) = \frac{1}{2}(m^2 +m)+m.
\]
We remark that this matroid argument does not apply to {\bf Case 1}, where $\bA/\BE^1$ may not be unimodular so $|\bA/ \BE^1| \le \cc(1,m-1)$ cannot be used in the last inequality.
\hfill\halmos 

\medskip

\noindent{\bf Tight examples.}
Proposition~\ref{prop:lowerbounds} provides an example of a bimodular matrix $\bA$ with a non-primitive column that satisfies $|\bA| = \cc(2,m)$; therefore, the upper bound in {\bf Case 2} is tight.
In the following paragraphs, we discuss the tightness of the bounds \eqref{eqe1e2}, \eqref{eqe1e2e3}, and \eqref{eqe1e2e3e4}, which consider when $\bA$ only contains primitive columns. 
We highlight two special cases: when $m=3$ and $|\bB^*|=2$ and when $m=5$ and $|\bB^*|=4$.
These cases are special because according to \eqref{eqe1e2}, \eqref{eqe1e2e3}, \eqref{eqe1e2e3e4} and {\bf Subcase 1.2}, they are the only cases when $|\bA|$ may equal $\cc(2,m)$.

When $|\bB^*|=2$, the bound~\eqref{eqe1e2} is attainable.
If $m\ge 4$ or $m=2$, then a tight example comes from deleting the non-primitive column from the example for $\cc(2,m)$ in Section \ref{seclowerbounds}.
This example is the vertex-edge incidence matrix of the directed complete graph on $m$ vertices together with the identity matrix and $m-1$ extra columns $\BE^1+\BE^i$ for $i=2,\dots, m$. 
In this example, $\bB^*$ corresponds to $[\BE^1-\BE^2|\BE^1+\BE^2]$.
For the special case when $m=3$, a tight example with $\cc(2,3) = 9$ columns is the matrix in Claim~\ref{claim_m=3}:
\[
\left[
\begin{array}{
c@{\hskip.1 cm}|@{\hskip.1 cm}
c@{\hskip.1 cm}|@{\hskip.1 cm}
c@{\hskip.1 cm}|@{\hskip.1 cm}
c@{\hskip.1 cm}|@{\hskip.1 cm}
c@{\hskip.1 cm}|@{\hskip.1 cm}
c@{\hskip.1 cm}|@{\hskip.1 cm}
c@{\hskip.1 cm}|@{\hskip.1 cm}
c@{\hskip.1 cm}|@{\hskip.1 cm}
c
}
1& 1 & 0 & 1 &   0 &   1 &   1 & 2 & 1 \\
0& 2 & 1 & 1 &   0 &   0 &   2 & 2 & 1\\
0& 0 & 0 & 0 &   1 &   1 &   1 & 1 & 1
\end{array}
\right]
\sim
\left[
\begin{array}{
r@{\hskip.1 cm}|@{\hskip.1 cm}
c@{\hskip.1 cm}|@{\hskip.1 cm}
c@{\hskip.1 cm}|@{\hskip.1 cm}
c@{\hskip.1 cm}|@{\hskip.1 cm}
r@{\hskip.1 cm}|@{\hskip.1 cm}
r@{\hskip.1 cm}|@{\hskip.1 cm}
c@{\hskip.1 cm}|@{\hskip.1 cm}
c@{\hskip.1 cm}|@{\hskip.1 cm}
c
}
 1& 1 & 0 & 1 &  -1 &   0 &   0 & 1 & 0 \\
-1& 1 & 1 & 0 &   0 &  -1 &   1 & 0 & 0\\
 0& 0 & 0 & 0 &   1 &   1 &   1 & 1 & 1
\end{array}
\right].
\]
This example is the vertex-edge incidence matrix of the directed complete graph on three vertices appended to an identity matrix and three extra columns $\BE^1+\BE^2,\BE^1+\BE^3,\BE^2+\BE^3$.

When $|\bB^*|=3$, the upper bound \eqref{eqe1e2e3} for $m\ge 5$ or $m=3$ can be achieved by the example of the vertex-edge incidence matrix of the directed complete graph with the identity matrix and $m-2$ extra columns $\BE^1+\BE^2-\BE^i$ for $i=3,\dotsc,m$. 
For $m=4$, the upper bound \eqref{eqe1e2e3} can be achieved by the previous example with one extra column $\BE^1+\BE^2-\BE^3-\BE^4$.
In these examples, one choice of $\bB^*$ is $[\BE^3|\BE^1-\BE^2|\BE^1+\BE^2-\BE^3]$.

When $|\bB^*|=4$, the upper bound \eqref{eqe1e2e3e4} is tight for $m=4,5$.
By setting $\bA$ to the right hand side~\eqref{BSpan}, we obtain a tight example when $m = 4$; here, $\bB^*=[\BE^1|\BE^2|\BE^3|\BE^1+\BE^2+\BE^3+2\BE^4]$.
For the special case when $m=5$, a tight example with $20=\cc(2,5)$ many columns consists of the twelve columns in the previous example for $m=4$ and eight extra columns:
%
\[
\left[
\begin{array}{
c@{\hskip.1 cm}|@{\hskip.1 cm}
c@{\hskip.1 cm}|@{\hskip.1 cm}
c@{\hskip.1 cm}|@{\hskip.1 cm}
c@{\hskip.1 cm}|@{\hskip.1 cm}
c@{\hskip.1 cm}|@{\hskip.1 cm}
c@{\hskip.1 cm}|@{\hskip.1 cm}
c@{\hskip.1 cm}|@{\hskip.1 cm}
c@{\hskip.1 cm}|@{\hskip.1 cm}
c@{\hskip.1 cm}|@{\hskip.1 cm}
c@{\hskip.1 cm}|@{\hskip.1 cm}
c@{\hskip.1 cm}|@{\hskip.1 cm}
c@{\hskip.1 cm}|@{\hskip.1 cm}
c@{\hskip.1 cm}|@{\hskip.1 cm}
c@{\hskip.1 cm}|@{\hskip.1 cm}
c@{\hskip.1 cm}|@{\hskip.1 cm}
c@{\hskip.1 cm}|@{\hskip.1 cm}
r@{\hskip.1 cm}|@{\hskip.1 cm}
r@{\hskip.1 cm}|@{\hskip.1 cm}
c@{\hskip.1 cm}|@{\hskip.1 cm}
c
}
1& 0 & 0 & 1 & 0 & 1 & 0 & 1 & 0 & 1 & 0 & 1 & 0& 1 & 0 & 1 &   0 &   0 &   1 & 1 \\
0& 1 & 0 & 1 & 0 & 0 & 1 & 1 & 0 & 0 & 1 & 1 & 0& 0 & 1 & 1 &   0 &   0 &   1 & 1 \\
0& 0 & 1 & 1 & 0 & 0 & 0 & 0 & 1 & 1 & 1 & 1 & 0& 0 & 0 & 0 &  -1 &   0 &   0 & 1 \\
0& 0 & 0 & 2 & 1 & 1 & 1 & 1 & 1 & 1 & 1 & 1 & 0& 0 & 0 & 0 &  -1 &  -1 &   1 & 1 \\
0& 0 & 0 & 0 & 0 & 0 & 0 & 0 & 0 & 0 & 0 & 0 & 1& 1 & 1 & 1 &   1 &   1 &   1 & 1
\end{array}
\right].
\]
This example is equivalent to the vertex-edge incidence matrix of the directed complete graph minus the column $\BE^1-\BE^2$, along with the identity matrix and six extra columns: $\BE^1+\BE^2-\BE^i$ for $i=3,4,5$ and  $\BE^1+\BE^2-\BE^3-\BE^4-\BE^5+\BE^i$ for $i=3,4,5$.
Here, $\bB^*$ corresponds to $[\BE^3|\BE^1-\BE^4|\BE^2-\BE^4|\BE^1+\BE^2-\BE^3]$.

The bound \eqref{eqe1e2e3e4} for $m\ge 6$ already shows that the maximal example in this case has at most $\cc(2,m)-1$ differing columns.
However, we believe that this bound can be improved. 
%
%
%
Our current maximal example found for $m=6$ has $\sfrac12\cdot(m^2+m)+m-2=\cc(2,m)-2=25$ columns.
The example is the vertex-edge incidence matrix of the directed complete graph with identity matrix except $\BE^1,\BE^2$, along with $-\BE^1+\BE^i+\BE^j$ for $i\ne j\in\{3,4,5,6\}$. 
Here, one choice of $\bB^*$ is $[\BE^3-\BE^4|\BE^5-\BE^6|-\BE^1+\BE^3+\BE^4|-\BE^1+\BE^5+\BE^6]$.
Our current maximal example found for $m\ge 7$ has $\sfrac{1}{2}\cdot(m^2+m)+m-3=\cc(2,m)-3$ columns.
The example is the vertex-edge incidence matrix of the directed complete graph without the column $\BE^1-\BE^2$, along with the identity matrix and the columns $\BE^1+\BE^2-\BE^i$ for $i=3,\dotsc,m$.
Here, $\bB^*$ corresponds to $[\BE^3|\BE^1-\BE^4|\BE^2-\BE^4|\BE^1+\BE^2-\BE^3]$.
%

\section{A proof of Proposition~\ref{prop:m2}.}\label{secm=2}

If $m = 1$, then $[1|\cdots|\Delta]$ is the unique maximal $\Delta$-modular matrix with differing columns (up to multiplication by $-1$).
Thus, the result holds for $m =1$.

Suppose $m=2$.
Let $\bA \in \Z^{2 \times n}$ be a $\Delta$-modular matrix with differing columns that satisfies $\rank \bA=2$ and $|\bA| = \cc(\Delta,2)$. 
Let $[\bb^1|\bb^2] \subseteq \bA$ satisfy $|\det [\bb^1|\bb^2]| = \Delta$.
Each column $\ba \in \bA$ can be written as $ \ba = v_1\bb^1+v_2\bb^2$ for $v_1, v_2 \in [-1,1]$. 
Otherwise, say if $|v_1| > 1$, then we derive the contradiction $|\det [\ba |\bb^2]| = |v_1| \cdot |\det [\bb^1|\bb^2]| > \Delta$.
After possibly multiplying columns of $\bA$ by $-1$, we assume that $v_2 \in [0,1]$ for each column $v_1\bb^1+v_2\bb^2 \in \bA$.

Set $\Pi := \{v_1 \bb^1+v_2 \bb^2 \in \Z^2 \colon v_1, v_2 \in [0,1)\}$.
It is well known that $|\Pi| = |\det[\bb^1|\bb^2]| = \Delta$; see~\cite[\S VII]{barv2002}.
Partition $\Pi$ as $\Pi = \{\bzero\} \cup \Pi^1 \cup \Pi^2 \cup \Pi^{{\rm int}}$, where
\[
\begin{array}{lcllll}
\Pi^1 &:=&  \left\{v_1 \bb^1  \in \Pi \colon v_1 \in (0,1)\right\},
&\\[.1 cm]
\Pi^{2} & := & \left\{ v_2 \bb^2 \in \Pi \colon v_2 \in (0,1)\right\},&\\[.1 cm]
\Pi^{{\rm int}} & := & \left\{v_1 \bb^1  + v_2 \bb^2 \in \Pi \colon (v_1,v_2) \in (0,1)^2 \right\} &.
\end{array}
\]

After multiplying columns by $-1$, we assume that if $v_1\bb^1 \in A$ for some $v_1 \in [-1,1]$, then $v_1 \ge 0$.
Hence, we assume $\bA \cap (\Pi^1+\{\bzero, - \bb^1\}) \subseteq  \Pi^1 $.\footnote{For sets $X, Y \subseteq \R^d$, the {\it Minkowski sum} is $X + Y := \{\bx + \by \colon \bx \in X, ~\by \in Y\}$. For $\by \in \R^d$, we write $X+\by$ instead of $X  +\{\by\}$.}
We partition $\bA \setminus \Pi$ as follows:
\begin{align}
\bA\setminus \Pi = &\ \left(\bA \cap \left(\Pi^{{\rm int}} -\bb^1\right)\right)\notag\\[.1 cm]
\cup &\ \left(\bA \cap \left(\Pi^{1}+\{ \bb^2, -\bb^1+\bb^2\}\right)\right) \cup \left(\bA \cap \{\bb^2, \bb^1+\bb^2, -\bb^1+\bb^2\}\right)\label{m2b1} \\[.1 cm]
\cup &\ \left(\bA \cap \left(\Pi^{2}+\{ \pm \bb^1\}\right)\right) \cup \{\bb^1\}.\label{m2b2}
\end{align}

Suppose $v_1\bb^1+\bb^2\in \bA$ for some $v_1 \in [-1,1]$; it follows that if $w_1\bb^1+\bb^2\in \bA$ for some $w_1 \in [-1,1]$, then $|v_1 - w_1| \le 1$. 
Indeed, otherwise we obtain the contradiction $|\det[v_1\bb^1+\bb^2|w_1\bb^1+\bb^2] | = |v_1-w_1|\cdot|\det[\bb^1|\bb^2]| > \Delta$.
This implies that the set in~\eqref{m2b1} has cardinality at most $ |\Pi^1| + 2$; furthermore, the cardinality is equal to $|\Pi^1|+2$ if and only if $\bA$ contains two columns $v_1\bb^1+\bb^2$ and $(v_1+1)\bb^1+\bb^2$ for some $v_1 \in [-1,1]$.
Similarly, it can be shown that the set in~\eqref{m2b2} has cardinality at most $ |\Pi^2| + 2$. 

We claim that either the set in~\eqref{m2b1} has cardinality at most $ |\Pi^1| + 1$ or the set in~\eqref{m2b2} has cardinality at most $ |\Pi^2| + 1$.
Assume that the set in~\eqref{m2b1} has cardinality $ |\Pi^1| + 2$.
Thus, $\bA$ contains two columns of the form $v_1 \bb^1+\bb^2$ and $(v_1+1)\bb^1+\bb^2$. 
Replace $\bb^2$ with $v_1\bb^1+\bb^2$; the result is another basis in $\bA$ of absolute determinant $\Delta$.
After this replacement, we can assume that $\bA$ contains the columns $\bb^2$ and $\bb^1+\bb^2$. 
The matrix $\bA$ cannot contain a column of the form $-\bb^1+v_2 \bb^2$ for $v_2 > 0$ otherwise $|\det[\bb^1+\bb^2|-\bb^1+v_2\bb^2]| = |1+v_2| \cdot|\det[\bb^1|\bb^2]| > \Delta$.
Hence, the set in~\eqref{m2b2} is contained in $\left(\bA \cap (\Pi^2 + \bb^1)\right)\cup \{\bb^1\}$, which contains at most $|\Pi^2|+1$ many elements. 
In other words, the union of the sets in~\eqref{m2b1} and~\eqref{m2b2} has cardinality at most $|\Pi^1|+|\Pi^2|+3$.

The column $\bzero $ is in $\Pi \setminus \bA$, and $|\Pi| = \Delta$.
Therefore, $|A \cap \Pi| = |A \cap (\Pi^{{\rm int}} \cup \Pi^1 \cup \Pi^2)| \le |\Pi| - 1$. 
The set $\bA \cap \left(\Pi^{{\rm int}} -\bb^1\right)$ is contained in a translation of $\Pi^{{\rm int}}$, so $|\bA \cap \left(\Pi^{{\rm int}} -\bb^1\right)| \le |\Pi^{{\rm int}}|$. 
By combining our upper bounds and applying the claim in the previous paragraph, we obtain
\[
|\bA| = |\bA \cap \Pi| +  |\bA \setminus \Pi| \le (|\Pi|-1) + (|\Pi^{{\rm int}}|+|\Pi^1|+|\Pi^2|+3) = 2(|\Pi|-1)+3 = 2\Delta+1.
\]
The equation $\cc(\Delta,2) = 2\Delta+1$ then follows from Proposition~\ref{prop:lowerbounds}.

\section{A proof of Theorem~\ref{thmLuze}.}\label{secGen}

Let $\bA \in \Z^{m\times n}$ be a $\Delta$-modular matrix with $\rank \bA=m$ and differing columns that satisfies $\cc(\Delta,m) = |\bA|$.
Let $\bB \subseteq \bA$ satisfy $|\det \bB| = \Delta$.
If $\Delta \le 2$, then Theorem~\ref{thmLuze} follows from Theorem~\ref{thm:bimodular} and Heller's result.
Therefore, assume $\Delta \ge 3$.
We can assume that $\bB$ is in Hermite Normal Form (See~\cite[\S 4.1]{AS1986}):\footnote{We use the convention that blank entries in a matrix are zero.}
\[
\bB = 
\left[
\begin{array}{c@{\hskip.1cm}|@{\hskip.1cm}c@{\hskip.1cm}|@{\hskip.1cm}c}
\bb^1 & \cdots & \bb^m
\end{array}
\right] \sim 
\left[
\begin{array}{c@{\hskip .1 cm}|@{\hskip .1 cm}c@{\hskip .1 cm}|@{\hskip .1 cm}c@{\hskip .1 cm}|@{\hskip .1 cm}c}
\mathbb{I}_{m-k} &  * & \cdots & *  \\
& \delta_1 &  \ddots & \vdots \\
& & \ddots &  * \\
&&&\delta_k
\end{array}
\right],
\]
where $\delta_1, \dotsc, \delta_k \ge 2$, $\prod_{i=1}^k \delta_i = \Delta$, and for each $\bb^i = (b^i_1, \dotsc, b^i_m)$, we have $0 \le b^i_j < b^i_i$ for all $j = 1, \dotsc, i-1$ and $ b^i_j =0$ for all $j=i+1,\dotsc, m$. 

Each column $\ba = (a_1, \dotsc, a_m) \in  \bA \setminus \bB$ satisfies $|a_m| \le \delta_k$ because $|a_m|\cdot \prod_{i=1}^{k-1} \delta_i = |\det[\bb^1|\cdots|\bb^{m-1}|\ba]| \le \Delta$. 
After possibly multiplying columns by $-1$, we can assume that $a_m \in \{0, \dotsc, \delta_k\}$ for all columns $\ba$.
For $r \in \Z$, define
\[
\bA[r]:= \big\{ (a_1, \dotsc, a_m) \in \bA \colon a_m = r\big\}.
\]
For each prime number $p$, define
\[
\hat{\bA}[p]:=\bigcup_{\substack{i=1, ~ p\mid i}}^{\delta_k} \bA[i].
\]
%
By applying a union bound, we see that
\begin{equation}\label{eqGenInduction}
\cc(\Delta,m)  = |\bA| \le |\bA[1]| + \sum_{\substack{p=2,\\p~\text{prime}}}^{\delta_k} ~\left|\hat{\bA}[p]\right|.
\end{equation}
We use~\eqref{eqGenInduction} to upper bound $|\bA|$ in terms of $\cc(1,m), \dotsc,$ and $\cc(\Delta-1,m)$.
Our analysis distinguishes between the cases $k =1$ and $k \geq 2$.

\smallskip
\noindent{\bf Case 1.}
Assume that $k = 1$. 
Recall $\Delta \ge 3$. 
For this range of $\Delta$, Glanzer et al.~\cite[Subsection 3.2]{GWZ2018} showed that
    \begin{align}\label{eqn:bound4}
        \cc(\Delta,m) &\le \sum_{\substack{p=2,\\ p~\text{prime}}}^{\Delta} \cc\left(\left\lfloor\frac{\Delta}{p}\right\rfloor,m\right) + 2\cc\left(\left\lfloor\frac{\Delta}{2}\right\rfloor,m\right).
    \end{align}
\smallskip

\noindent{\bf Case 2.} Assume that $k \ge 2$.
For a prime $p \in \{2, \dotsc, \delta_k\}$ and an integer $i$ divisible by $p$, we can divide the $m$th row of $\bA[i]$ by $p$. 
We have
\begin{equation}\label{eqGenBoundNot1}
\left|\hat{\bA}[p]\right|  ~\le~ \left|\left[
\begin{array}{c@{\hskip.1cm}|@{\hskip.1cm}c@{\hskip.1cm}|@{\hskip.1cm}c@{\hskip.1cm}|@{\hskip.1cm}c}
\bb^1& \cdots&\bb^{m-1}& \hat{\bA}[p]
\end{array}\right]\right| ~\le~ \cc\left(\left\lfloor\frac{\Delta}{p}\right\rfloor,m\right)
\end{equation}
because the columns in the middle expression form a $\floor{\sfrac{\Delta}{p}}$-modular matrix with rank-$m$.

Consider $\bA[1]$.
For each $\ba = (a_1, \dotsc, a_m) \in A[1]$, we have
\begin{align*}
\big|
\det [\bb^1 ~|\cdots| ~\bb^{m-2}~ |~\ba~|~\bb^m] \big|
= &\
%
\left|
\det
\left[
\begin{array}{c@{\hskip .1 cm}|@{\hskip .1 cm}c@{\hskip .1 cm}|@{\hskip .1 cm}c@{\hskip .1 cm}|@{\hskip .1 cm}c@{\hskip .1 cm}|@{\hskip .1 cm}c@{\hskip .1 cm}|@{\hskip .1 cm}c}
\mathbb{I}_{m-k} & * &\cdots& *& *& *  \\ 
  & \delta_1 & \cdots & * & * & *\\
  &  & \ddots & \vdots & \vdots& \vdots\\
& & &\delta_{k-2} & * & *\\
& &  & & a_{m-1} & b^m_{m-1}  \\
&&&&1 & \delta_k\\
\end{array}
\right]
\right|\\[.25 cm]
=&\  \big|a_{m-1}\delta_k - b^m_{m-1}\big|\prod_{i=1}^{k-2} \delta_i\\
\le&\ \Delta.
\end{align*}
Thus, $|a_{m-1}\delta_k - b^m_{m-1}| \le \delta_{k-1}\delta_k$.
The Hermite Normal Form assumption implies $0 \le b^m_{m-1} < \delta_k$, so we have $|a_{m-1}| \le \delta_{k-1}$.
For each $r \in \Z$, define the column set 
\[
\bA[1,r] := \big\{ (a_1, \dotsc, a_m) \in \bA[1] \colon |a_{m-1}| = r\big\}.
\]
With these sets, we can upper bound $|A[1]|$:
\[
|\bA[1]| \le |\bA[1,\delta_{k-1}]| + \sum_{i=0}^{\floor{\log_2(\delta_{k-1}-1)}} \bigg|\bigcup_{\substack{s\in \Z,~s~\text{odd}}} \bA[1, s2^i]\bigg|.
\]
The matrix $[\bb^1| \cdots| \bb^{m-1}| \bA[1, \delta_{k-1}]]$ is $\Delta$-modular of full row rank, and the $(m-1)$st row is divisible by $\delta_{k-1}$; dividing the $(m-1)$st row by $\delta_{k-1}$ shows $|\bA[1, \delta_{k-1}]| \le \cc(\sfrac{\Delta}{\delta_{k-1}},m)$.

For each $i = 0, \dotsc, \floor{\log_2(\delta_{k-1}-1)}$, define
\[
\overline{\bA}[i] := \bigcup_{\substack{s\in \Z,~s~\text{odd}}} \bA[1, s2^i].
\]
If $\overline{\bA}[i] = \bA[1, s2^i]$ for a single odd integer $s$, then perform the following elementary operation to the full row rank matrix $\left[\bb^1|\cdots|\bb^{m-1}|\overline{\bA}[i]\right]$: subtract $s2^i$ times the $m$th row from the $(m-1)$st row.
The $(m-1)$st row of the resulting matrix is divisible by $\delta_{k-1}$, so $|\overline{\bA}[i]| \le \cc(\sfrac{\Delta}{\delta_{k-1}}, m)$.
Suppose $\overline{\bA}[i] \neq \bA[1, s2^i]$ for a single odd integer $s$.
For each column $\ba = (a_1, \dotsc, a_{m-2}, s2^i,1)\in\overline{\bA}[i]$, if we add $2^i$ times the $m$th row to the $(m-1)$st row, then the resulting column is $(a_1, \dotsc, a_{m-2}, (s+1)2^i,1)$; in particular, the $(m-1)$st entry is divisible by $2^{i+1}$ because $s$ is odd. 
Perform this elementary operation to the full rank matrix $\left[\bb^1| \cdots|\bb^{m-2}| \overline{\bA}[i]\right]$: add $2^i$ times the $m$th row to the $(m-1)$st row.
The $(m-1)$st row of the resulting matrix is divisible by $2^{i+1}$.
Hence, $|\overline{\bA}[i]|\le \cc(\floor{\sfrac{\Delta}{2^{i+1}}},m)$.

Substituting~\eqref{eqGenBoundNot1} and our bounds for $|\bA[1, \delta_{k-1}]|$ and $\left|\overline{\bA}[i]\right|$ into~\eqref{eqGenInduction}, we see that
    \[
    \cc(\Delta,m) \le \sum_{\substack{p=2,\\ p~\text{prime}}}^{\delta_k} \cc\left(\left\lfloor\frac{\Delta}{p}\right\rfloor,m\right) + \cc\left(\frac{\Delta}{\delta_{k-1}},m\right) + \sum_{\ell=0}^{\floor{\log_2(\delta_{k-1}-1)}} \hspace{-.25 cm}\max \left\{\cc\left(\frac{\Delta}{\delta_{k-1}},m\right), \cc\left(\left\lfloor\frac{\Delta}{2^{\ell+1}}\right\rfloor,m\right)\right\}.
    \]
    If $\delta_{k-1}=2$, then $\floor{\log_2(\delta_{k-1}-1)}=0$ and
    \begin{equation}\label{eqn:bound1}
        \cc(\Delta,m) \le \sum_{\substack{p=2,\\ p~\text{prime}}}^{\delta_k} \cc\left(\left\lfloor\frac{\Delta}{p}\right\rfloor,m\right) + 2\cc\left(\frac{\Delta}{2},m\right).
    \end{equation}
    If $\delta_{k-1}=3$, then $\floor{\log_2(\delta_{k-1}-1)}=1$ and
    \begin{equation}\label{eqn:bound2}
       \cc(\Delta,m) \le \sum_{\substack{p=2,\\ p~\text{prime}}}^{\delta_k} \cc\left(\left\lfloor\frac{\Delta}{p}\right\rfloor,m\right) + 2\cc\left(\frac{\Delta}{3},m\right) + \cc\left(\left\lfloor\frac{\Delta}{2}\right\rfloor,m\right).
    \end{equation}
    If $\delta_{k-1}\ge 4$, then
    \begin{equation}\label{eqn:bound3}
        \cc(\Delta,m) \le \sum_{\substack{p=2,\\ p~\text{prime}}}^{\delta_k} \cc\left(\left\lfloor\frac{\Delta}{p}\right\rfloor,m\right) + 2\cc\left(\frac{\Delta}{\delta_{k-1}},m\right) + \sum_{\ell=0}^{\floor{\log_2(\delta_{k-1}-1)}-1} \cc\left(\left\lfloor\frac{\Delta}{2^{\ell+1}}\right\rfloor,m\right).
    \end{equation}
    This completes {\bf Case 2.}

We now use the two cases to bound $\cc(\Delta, m)$. 
Define $\cg(\Delta,m):=\sfrac{1}{2}\cdot(m^2+m)\Delta^2$. 
We use induction on $\Delta$ to show $\cc(\Delta, m)\le \cg(\Delta, m)$.
Glanzer et al. demonstrated that $\cc(\Delta,m)\le \cg(\Delta,m)$ for $\Delta \le 3$; see~\eqref{eqOldBounds}.
Thus, we assume that $\Delta \ge 4$ and $\cc(\delta,m) \le \cg(\delta,m)$ for each $\delta < \Delta$.
To prove $\cc(\Delta, m)\le \cg(\Delta, m)$, it suffices to upper bound the right hand sides of~\eqref{eqn:bound4},~\eqref{eqn:bound2}, and~\eqref{eqn:bound3} by $\cg(\Delta, m)$.
Given that $\Delta \ge 4$, we do not need to bound~\eqref{eqn:bound1} because it is less than~\eqref{eqn:bound4}.

Using the definition of $\cg(\delta, m)$ and the prime zeta function $\cp(s):= \sum_{p\text{ prime}}\sfrac{1}{p^s}$, we arrive at the bound
\begin{equation}\label{eqprimeZeta}
\sum_{\substack{p=2,\\ p~\text{prime}}}^{\Delta} \cg\left(\left\lfloor\frac{\Delta}{p}\right\rfloor,m\right)\le \sum_{\substack{p=2,\\ p~\text{prime}}}^{\Delta} \frac{1}{p^2} \cg(\Delta,m) \le \cp(2) \cg(\Delta,m).
\end{equation}
We know $\cp(2)<\sfrac{1}{2}$~\cite[A085548]{sloane2018line}.
We extend~\eqref{eqn:bound4} using the induction hypothesis,~\eqref{eqprimeZeta}, and the definition of $\cg(\floor{\Delta/2},m)$:
    \[
        \cc(\Delta,m) 
        %
        \le \left(\cp(2)+\frac{1}{2}\right)\cg(\Delta,m)
        <\cg(\Delta,m). 
    \]
Similarly, we extend~\eqref{eqn:bound2}:
\[
       \cc(\Delta,m) 
        \le \left(\cp(2) + \frac{2}{9} + \frac{1}{4}\right)\cg(\Delta,m)
        <\cg(\Delta,m).
\]
We extend~\eqref{eqn:bound3} using $\sum_{i=0}^{t-1} \sfrac{1}{4^{i+1}} = \sfrac{1}{3}\cdot \left(1-\sfrac{1}{4^t}\right)<\sfrac{1}{3}$, the induction hypothesis and~\eqref{eqprimeZeta}:
\[
%
        %
        \cc(\Delta, m) \le \left(\cp(2) + \frac{2}{\delta_{k-1}^2} + \frac{1}{3}\right)\cg(\Delta,m)
        \le \left(\cp(2) + \frac{1}{8} + \frac{1}{3}\right)  \cg(\Delta,m)
        <\cg(\Delta,m).
\]
%
%
\hfill \halmos

\medskip
As a final remark, the term $\Delta^2$ in Theorem~\ref{thmLuze} comes from our ability to bound $\cg(\Delta,m)$ via induction.
If we apply our proof analysis to some other upper bound $\overline{\cg}(\Delta,m) = m^2\Delta^q$, where $1\le q\le 2$, then it suffices for $q$ to belong to the following set:
\[
\left\{q \in \R_{\ge 0} \colon \max \left\{\frac{2}{3^q} + \frac{1}{2^q}, \frac{2}{4^q} + \frac{1}{2^q-1}, \frac{2}{2^q}\right\} \le 1 - \cp(q)\right\}.
\]
%
From numerical computation, $q=1.95$ is in this set.
Thus, $\cc(\Delta,m) \le m^2 \Delta^{1.95}$.
We chose to present $m^2\Delta^2$ for simplicity.
%

An open question is whether $\cc(\Delta,m) \le \mathfrak{h}(m)\Delta$ for a polynomial $\mathfrak{h}$.
Since posting the first version of this manuscript, we became aware of geometric arguments by Gennadiy Averkov and Matthias Schymura that $\cc(\Delta,m)\le \mathfrak{f}(m)\Delta$, where $\mathfrak{f}$ is super-polynomial~\cite{AS2021}.
In fact, we can give a quick argument showing that $\cc(\Delta,m) \le 3^m\Delta$.
If $\bB \subseteq \bA$ satisfies $|\det \bB| = \Delta$, then every $\ba \in \bA$ satisfies $\|\bB^{-1}\ba\|_{\infty} \le 1$.
Hence, $\bA \subseteq \Pi + B\{-1,0,1\}^m$, where $\Pi := \{\bB \bv \in \Z^m \colon \bv \in [0,1)^m\}$.
From this, we see that $|\bA| \le 3^m \cdot |\Pi| = 3^m\Delta$.
%

\section{A proof of Theorem~\ref{thmProx}.}\label{secProx}

Recall that when it comes to bounding $\pi$ the matrix $\bA$ does not necessarily have differing columns as is the case with the other results in this manuscript.

Let $\bx^*$ be a vertex of $\LP$ satisfying $\pi = \min_{\bz \in\IP} ~~\|\bx^* - \bz\|_1$.
By standard $\LP$ results,  there exists a vector $\bc \in \R^n$ such that $\bx^*$ is the unique maximizer of $\bx \to \bc^\intercal \bx$ over $\LP$; see~\cite[Chapter 3]{CCZ2014}.
By possibly perturbing $\bc$, we can assume that there is unique maximizer $\bz^*$ of $\bz \to \bc^\intercal \bz$ over $\IP$.
Note that $\bz^*$ is a vertex of $\conv \IP$.
Let $k  \le 2\cc(\Delta, m)+1$ denote the number of distinct columns of $\bA$.
By applying~\cite[Theorem 2]{PSW2020} with $T = \emptyset$ and $B = \bA$, there exists a matrix $\bW \in \Z^{k\times n}$ such that
\[
\conv \IP = \conv\{\bx \in \LP \colon \bW \bx \in \Z^k\}.
\]
The previous equation implies $\bz^*$ is also the unique maximizer of $\bz \to \bc^\intercal \bz$ over the \emph{mixed integer} linear set $\{\bx \in \LP \colon \bW \bx \in \Z^k\}$, which has $k$ many integer constraints.

Consider the difference vector $\bx^* - \bz^*$.
From the equation $\pi = \min_{\bz \in\IP} ~~\|\bx^* - \bz\|_1$, it follows that $\pi \le \|\bx^* - \bz^*\|_1$.
The difference vector was first analyzed by Cook et al.~\cite{CGST1986} and later in~\cite{GraProx,HochConv,LPSX2020,PWW2018,WerSep}.
The proof of mixed integer proximity in~\cite[Theorem 2]{PWW2018} established that
\[
\bx^* - \bz^* = \sum_{i=1}^n \lambda_i \bu^i,
\]
where $\bu^1, \dotsc, \bu^n \in \Z^n$ and $\lambda_1, \dotsc, \lambda_n \ge 0$ satisfy $\sum_{i=1}^n \lambda_i \le k$.
The result~\cite[Claim 8]{LPSX2020} demonstrated that $\|\bu^i\|_{1} \le (m+1)\Delta$ for each $i = 1, \dotsc, n$.
Therefore, $\pi \le \|\bx^* - \bz^*\|_1 \le (m+1)\Delta k$.
The result now follows from Theorem~\ref{thmLuze}.
\hfill\halmos
\section*{Acknowledgments}
J. Lee was supported in part by ONR grant N00014-21-1-2135 and AFOSR grant FA9550-19-1-0175. 
J. Paat was supported by a Natural Sciences and Engineering Research Council of Canada (NSERC) Discovery Grant [RGPIN-2021-02475].
The authors would like to thank Jim Geelen for making us aware of the results in~\cite{GNW2021}.
We thank Stefan Kuhlmann for his comments on our lower bound construction, and we thank Rico Zenklusen for his comments regarding~\cite{CWZ2018}.
%


\bibliographystyle{informs2014}
\bibliography{references}

\end{document}